\documentclass[11pt, oneside]{article}   	
\usepackage[margin = 1in]{geometry}                		
\usepackage{graphicx}				
\usepackage{mathrsfs}
\usepackage{float}
\usepackage{pillowmath}
\usepackage{pgfplots}
\usepackage{caption}
\usepackage{nicematrix}
\usepackage{tikz-qtree}
\usepackage{tikz}
\usepackage{tikz-cd}
\usepackage{authblk}
\usepackage{thm-restate}

\newtheorem{theorem}{Theorem}[section]
\newtheorem{lemma}[theorem]{Lemma}
\newtheorem{proposition}[theorem]{Proposition}

\newtheorem{remark}[theorem]{Remark}
\newtheorem{definition}[theorem]{Definition}
\newtheorem{example}{Example}

\makeatletter
\newcommand{\MSC}[1]{%
  \begingroup
    \renewcommand\thefootnote{}
    \footnote{\textbf{2020 Mathematics Subject Classification.} #1}%
    \addtocounter{footnote}{-1}
  \endgroup
}
\newcommand{\keywords}[1]{%
  \begingroup
    \renewcommand\thefootnote{}%
    \footnote{\textbf{Key words and phrases.} #1}%
    \addtocounter{footnote}{-1}%
  \endgroup
}
\makeatother
  
  \ExplSyntaxOn
\NewExpandableDocumentCommand { \ValuePlusOne } { m } 
  { \int_eval:n { \int_use:c { c @ #1 } + 1 } }
\NewExpandableDocumentCommand { \Sec } { m } 
  { \fp_eval:n { secd ( #1 ) } }
\NewDocumentCommand { \Rot } { m }
  { 
    \hbox_to_wd:nn { 1 em }
      { 
        \hbox_overlap_right:n 
          { 
            \skip_horizontal:n { \fp_to_dim:n { 7 * cosd (45) } } 
            \rotatebox{45}{#1}
          } 
      } 
  }
\ExplSyntaxOff

\NewDocumentCommand { \MixedRule } { m }
  {
    \begin{tikzpicture}
    \coordinate (a) at (2-|#1) ;
    \coordinate (b) at (1-|#1) ;
    \draw (a) -- ($(a)!\Sec{90-45}!45-90:(b)$) ;
    \draw (2-|#1) -- (\ValuePlusOne{iRow}-|#1) ;
    \end{tikzpicture}
  }
  
   \newtheorem{innercustomthm}{Theorem}
\newenvironment{customthm}[1]
  {\renewcommand\theinnercustomthm{#1}\innercustomthm}
  {\endinnercustomthm}
  
\newtheorem{innercustomprop}{Proposition}
\newenvironment{customprop}[1]
  {\renewcommand\theinnercustomprop{#1}\innercustomprop}
  {\endinnercustomprop}

\newtheorem{innercustomlem}{Lemma}
\newenvironment{customlem}[1]
  {\renewcommand\theinnercustomlem{#1}\innercustomlem}
  {\endinnercustomlem}
  
\newtheorem{innercustomex}{Example}
\newenvironment{customex}[1]
  {\renewcommand\theinnercustomex{#1}\innercustomex}
  {\endinnercustomex}

\allowdisplaybreaks

\usepackage{hyperref}
\hypersetup{
    colorlinks=true,
    linkcolor=blue,
    filecolor=magenta,      
    urlcolor=cyan,
    citecolor=blue
}

\usepackage{thmtools}
\usepackage{thm-restate}

\title{How should we aggregate ratings? Accounting for personal rating scales via Wasserstein barycenters}
\date{}							

\begin{document}

\author{Daniel Raban\footnote{Email: danielraban@berkeley.edu}}
\affil{Department of Statistics, University of California, Berkeley}

\maketitle



\abstract{\footnotesize A common method of comparing items is to collect numerical ratings on a linear scale and compare the average rating for each item. However, averaging ratings does not account for people rating according to differing personal rating scales. With this in mind, we investigate the problem of calculating aggregate numerical ratings from individual numerical ratings and propose a new, non-parametric model for the problem. We show that, with minimal modeling assumptions, the standard average is inconsistent for estimating the quality of items. Analyzing the problem of heterogeneous personal rating scales from the perspective of optimal transport, we derive an alternative rating estimator, which we show is asymptotically consistent almost surely and in $L^p$ for estimating quality, with an optimal rate of convergence. Further, we generalize Kendall's W, a non-parametric coefficient of preference concordance between raters, from the special case of rankings to the more general case of arbitrary numerical ratings. Along the way, we prove Glivenko--Cantelli-type theorems for uniform convergence of the cumulative distribution functions and quantile functions for Wasserstein-2 barycenters on $[0,1]$.}

\MSC{62G05 (Primary); 49Q22, 62G20, 62H20, 62F07 (Secondary)}


\section{Introduction}

A common method of comparing items is to collect numerical ratings on a range of, say, 1 to 10, and to then average the ratings. Usage of numerical ratings is almost a hundred years old \cite{thurstone1928attitudes, likert1932technique}, and analysis of ratings data is ubiquitous, appearing in psychology, the study of consumer preferences, natural language processing, and more.


Unfortunately, the actual numerical values in ratings data are essentially meaningless if taken at face value. For instance, the
difference in quality between items rated 2 and 4 may not be the same as the difference in quality between items rated 4 and 6. Averaging the ratings for an item forces this interpretation on us, which, in general, may not be appropriate. \cite{stark2014evaluation} raises this criticism for student evaluations of teaching, and in the same context, \cite{mccullough2011analysing} argues that proportions are more appropriate than averages of rating data.

The central issue we will consider in this work---which we address in an attempt to alleviate the above concerns---is that people's ratings do not only provide noisy versions of some unknown ``quality score'' for an item, for which taking a standard average would be appropriate. Rather, people provide ratings on different \emph{scales}. In other words, each rater has their own, personal distribution of ratings that they follow, so the values from these different distributions are not directly comparable.

\begin{example}
Consider two users of a movie rating website, one of whom usually rates movies between 3 to 7 and the other of whom only rates half of the movies they watch 9 and the other half 10. If the two users rate movie A as 9, then their ratings do not mean the same thing. The first user considers movie A to be better in quality than the vast majority of movies, while the second user is saying that movie A is in the bottom half of movies.
\end{example}

In situations such as the above example, where we have access to \emph{multiple ratings by each user}, we can understand users' personal rating scales by looking at their personal distribution of ratings across multiple items. 
If there are $M$ items and $n$ users, then (in the ideal situation where all users have rated all items) the statistician has access to an entire matrix $(r_{i,j})_{1 \leq i \leq n, 1 \leq j \leq M}$ of ratings.
\begin{center}
\begin{NiceTabular}{ccccc}
 & \Rot{Inception} & \Rot{Oppenheimer} & \Rot{Tenet} & \Rot{The Prestige} \\ \hline
Alice & 5 & 6 & 4 & 8 \\ \hline  
Bob & 9 & 8 & 10 & 8 \\ \hline  
\CodeAfter
\MixedRule{2}
\MixedRule{3}
\MixedRule{4}
\MixedRule{5}
\MixedRule{6}
\end{NiceTabular}
\end{center}
Thus, the statistician has access to an empirical distribution $\mu_j$ of ratings for each user, given by the rows of such a matrix. Averaging the values of each column of this matrix, i.e., the average of ratings for each item, neglects the wealth of information available. A clear understanding of each user's personal rating distribution allows us to compare users' ratings across differing rating scales. For the rest of this paper, we will continue to use the terminology associated to ratings websites, i.e., ``users'' rating ``items.''
\paragraph*{Our contributions.} The purpose of this paper is to establish a framework for and attempt to address the issue of users rating on personal scales. The first novel contribution of this paper is to provide a new, non-parametric framework for the ratings problem which incorporates the notion of personal rating scales. Rather than modeling how ratings are generated, this framework instead specifies properties defining population consensus.

The second novel contribution of this paper to provide a better-informed scheme for calculating aggregate ratings for each item by using users' empirical rating distributions. In this vein, we will show that for our model, the standard average of ratings for each item is asymptotically \emph{inconsistent} for estimating the numerical quality scores of items. Furthermore, we will propose an alternative method of producing aggregate ratings, which we show is asymptotically consistent for estimating the numerical quality of items, with an optimal $1/\sqrt n$ rate of convergence.

The key insight in both of these aspects will be consideration of the ratings problem from the perspective
of optimal transport, which naturally arises due to the necessity of converting ratings from individualized scales to a common scale.


The third novel contribution of this paper is to propose a generalization of Kendall's W statistic. Kendall's W measures the degree of agreement in a collection of ranked preference lists. Our framework for studying ratings includes rankings as a special case and suggests a natural generalization to the case of numerical ratings, yielding a statistic that measures the degree of agreement among rating profiles.
\paragraph*{Related work.} 
Our work builds on the theory of Wasserstein barycenters, a notion of averaging distributions common in the optimal transport literature (e.g.\ \cite{rabin2012wasserstein, le2017existence, chewi2020gradient}). Standard texts on the classical theory of optimal transport include \cite{rachev1998mass, villani2009optimal, villani2021topics}, and \cite{panaretos2020invitation} is a recent book which focuses on modern statistical applications of Wasserstein barycenters.

Our new framework for ratings resembles models from functional data analysis (e.g.\ \cite{gasser1995searching, gervini2004self} and \cite{10.1214/15-AOS1387}, which utilizes optimal transport in the context of estimation of point processes), in which we view a number of warped perspectives (user ratings, in our case) of an average behavior (numerical quality scores). We are then tasked with recovering the average behavior from the warped perspectives. A key difference in our case will be that, because users can disagree on preferences, we must necessarily allow for \emph{nonincreasing} warpings. See \cite{10.3150/13-BEJ585} for a general discussion of the use of Wasserstein barycenters in this way.


A popular alternative to using ratings data is to collect pairwise comparisons between various items under consideration. This method of comparisons has been well-studied in the statistics literature, most notably via the Bradley-Terry-Luce model \cite{bradley1952rank}, which has attracted recent interest devoted to extracting minimax rates (e.g.\ \cite{hendrickx2020minimax}). We will not directly compare our methods to BTL, as they take different types of data as inputs.

There is a significant amount of mathematical psychology literature on the topic of combining preferences. The use of rating data originated in this field \cite{thurstone1928attitudes, likert1932technique} and has persisted over the years \cite{saal1980rating, uher2018quantitative, nesbit2018statistics}.
One major legacy of Amos Tversky is a collection of papers which axiomatize the theory of preferences and study them from this perspective (e.g.\ \cite{tversky1969intransitivity,tversky1988contingent,tversky1993context}).



A somewhat recent body of work on ratings came from the Netflix prize, a public challenge to accurately predict user ratings given incomplete ratings data (see \cite{bennett2007netflix} for an introduction and early developments; \cite{koren2009bellkor} and \cite{toscher2009bigchaos} discuss the winning algorithm). We emphasize that aggregate ratings produce vector-valued summaries rather than completing the rating matrix, so they are not intended to impute missing user ratings. Rating-prediction methods, are therefore \emph{complementary}, rather than competing with aggregation, and future work could apply our methods after imputation.

\paragraph*{Definition of the rating estimator.} Our proposed rating estimator is given by the following 2-step computation: Suppose users have personal rating distributions $\mu_1,\dots,\mu_n$ with associated cumulative distribution functions $F_1,\dots,F_n$ (and hence generalized inverses $F_1^{-1},\dots,F_n^{-1}$, where $F_k^{-1}(y) := \inf \{ x : F_k(x) \geq y \}$).

\begin{itemize}

\item[] \textbf{Step 1.} (Primitive ratings): Given an item with user ratings $r_1,\dots,r_n$, define the \emph{primitive aggregate rating} for that item as
$$R_0(r_1,\dots,r_n) := \frac{1}{n^2} \sum_{k=1}^n \sum_{j=1}^n F_j^{-1} \circ F_k(r_k).$$
Calculate the primitive rating for all items, and obtain the distribution $\nu$ of primitive ratings for all items.

\item[] \textbf{Step 2.} (Final rescaling): Adjust each primitive rating $r$ by reporting the aggregate rating for that item as
$$R(r) := \frac{1}{n} \sum_{\ell=1}^n F_\ell^{-1} \circ F_{\nu}(r).$$

\end{itemize}
Step 2, applied to a primitive rating for each item, yields an aggregate rating for each item.

The interpretation is that $R_0$, which can be expressed as
$$R_0(r_1,\dots,r_n) = \frac{1}{n} \sum_{k=1}^n F_{\wh \mu}^{-1} \circ F_k(r_k),$$
where $\wh \mu$ is the Wasserstein-2 barycenter of $\mu_1,\dots,\mu_n$, first puts the ratings on a common numerical scale and then averages them. The maps $F_{\wh \mu}^{-1} \circ F_k$ are the optimal transport maps transporting user $k$'s personal rating distribution $\mu_k$ to the common/consensus scale distribution $\wh \mu$; these maps preserve the order of user $k$'s ratings while adjusting the scale of their rating distribution. The choice of the Wasserstein-2 barycenter $\wh \mu$ as the common scale to which we convert all personal rating distributions is canonical because it is the choice which minimizes the total squared transport distance to all the $\mu_k$.

An estimator in this context needs to recover both the scale/distribution and the order of population consensus numerical quality scores. If the users have differing preferences, then the resulting distribution $\nu$ of primitive ratings will not be exactly $\wh \mu$, which we view as the consensus scale that our aggregate ratings should follow, so we rescale the result to
$$R = F_{\wh \mu}^{-1} \circ F_{\nu} \circ R_0$$
to make the aggregate ratings more easily interpretable. Again, this rescaling preserves the ordering of the ratings and is hence an optimal transport map. We will see that when all users agree on preferences (but possibly rate on differing scales), the rating estimator $R$ equals the average. So it may be viewed as extending the average to the general case in which users disagree on preferences.
\paragraph*{Outline.} In what follows, we will investigate the estimator $R$ and its efficacy as a notion of aggregate ratings.

In Section~\ref{model_section}, we make explicit our new framework for the rating problem.

In Section~\ref{wasserstein_background}, we recall the necessary optimal transport background.

In Section~\ref{rating_estimation_with_universal_preferences}, we consider the special case where all users agree on preferences (but rate on differing scales). In Section~\ref{primitive_estimator}, we derive the primitive rating estimator $R_0$, using the geometry of the Wasserstein space of rating distributions. Proposition~\ref{R0_equals_A_and_is_consistent} provides almost sure and $L^p$ asymptotic consistency results for estimators of individual items'  quality.

In Section~\ref{rating_estimation_with_misaligned_preferences}, we extend our analysis and proofs of consistency to the general case in which users rate on different scales with differing preferences. Example~\ref{scaling_and_reversing} shows that the standard average of ratings is in general inconsistent for estimating the quality of items, Theorem~\ref{ratings_are_consistent} provides almost sure and $L^p$ consistency for our proposed rating estimator, and Theorem~\ref{rating_rates_of_convergence} provides rates of convergence for the estimation error.

In Section~\ref{incomplete-data}, we address the situation where not all users have rated all items. Theorem~\ref{incomplete-consistency-thm} provides consistency for our rating estimator in this case, with an optimal rate of convergence.

In Section~\ref{ranking-section}, we discuss the combining of preferences from ranked orders, which is a special case of our problem. Using this perspective, we propose two new statistics for measuring the concordance of user ratings and of their rating scales. Our statistic for user ratings generalizes Kendall's W statistic for concordance of rankings. Proposition~\ref{concordance-formulas} and Proposition~\ref{concordance-bounds} give properties of these statistics.

In Section~\ref{application_to_real_data}, we apply these estimators to real-world and simulated data.

\section{A new framework for the ratings problem} \label{model_section}

Any attempt to interpret ratings and compare them on a numerical scale presupposes that the items have some \emph{inherent numerical quality}, and so we model the problem of calculating aggregate ratings as that of estimation of all of these inherent qualities in such a way that they can be meaningfully compared. As ratings data is subjective, these inherent qualities are determined by population consensus. Section~\ref{modeling_assumptions} outlines the details of our framework. Section~\ref{identifiability_assumptions} provides minimal, non-parametric assumptions on how the notion of consensus determines numerical quality scores; these assumptions will make the model identifiable.

\subsection{Setup of the ratings problem} \label{modeling_assumptions}

From here onwards, ratings will be assumed to be values in $[0,1]$. We take the space of consensus numerical quality scores for items (unknown to the statistician) to be $[0,1]$, as well. The distribution on $[0,1]$ of consensus numerical quality scores for all items (unknown to the statistician, as well) will be $\mu$. For example, if we have $M$ items of differing quality under consideration, then $\mu = \frac{1}{M} \sum_{j=1}^M \delta_{x_j}$, where $x_j$ is the quality score for item $j$.

We first assume that all $n$ users have rated every item and that no two items share the same consensus quality. Thus, for the item with numerical quality $q$, the ratings $\phi_1(q),\dots,\phi_n(q)$ are well-defined $\mu$-a.e., yielding the personal \emph{scale functions} $\phi_1,\dots,\phi_n : [0,1] \to [0,1]$. Each user's rating distribution $\mu_k := (\phi_k)_*$ is the push-forward of the consensus distribution $\mu$ under their scale function $\phi_k$. In other words, each user transforms the consensus quality space $([0,1],\mu)$ into their own rating space $([0,1],\mu_k)$ via $\phi_k$.
$$\begin{tikzcd}[column sep=tiny]
& & ([0,1],\mu) \arrow[lldd,"\phi_1",swap]  \arrow[ldd,"\phi_2"] \arrow[rdd,"\phi_{n-1}",swap] \arrow[rrdd,"\phi_{n}"]  \\
\\
([0,1],\mu_1) & ([0,1],\mu_2) &  \cdots & ([0,1],\mu_{n-1}) & ([0,1],\mu_{n})
\end{tikzcd}$$
We work with the quality scores directly because identifying items with their quality scores embeds the items into $[0,1]$, yielding a natural ordering and confining the analysis to functions $[0,1] \to [0,1]$. In Section~\ref{incomplete-data}, we will extend this viewpoint to incorporate incomplete data.

Intuitively, we wish to capture the notion that all the ratings of an item are interpretations of a value on a space of consensus numerical quality scores. Consequently, our broad goal is to reconstruct the original space of quality scores after viewing $n$ distorted copies of the space, one from each user. This modeling setup resembles models of estimating functions in functional data analysis (e.g. \cite{gasser1995searching,gervini2004self}; \cite{10.1214/15-AOS1387} adapts these ideas to estimation of point processes), in which we view $n$ warped copies of some function, where the warping is done by altering the amplitude of the function's values and by warping the space via an \emph{increasing} function. In our context, the warping functions are the personal scale functions $\phi_k$. Crucially, however, users may disagree on preferences between items (i.e.\ we can have $\phi_1(x) < \phi_1(y)$ and $\phi_2(x) > \phi_2(y)$), so our $\phi_k$ need not be increasing, in general. Our model, viewed as an extension of these types of functional data analysis models to the case of non-increasing warp functions, will allow us to prove general consistency results which apply to situations where we do not only care about a finite set of values. Consequently, our theorems will not in general assume $\mu$ to be an atomic measure. \emph{However, for the purposes of aggregating ratings, the reader can safely assume that $\mu =  \frac{1}{M} \sum_{j=1}^M \delta_{x_j}$.}

We will assume the statistician has access to two kinds of data:
\begin{enumerate}[label = (\alph*)]

\item For $\mu$-a.e.\ $x \in [0,1]$ (that is, for the unique item with numerical quality score $x$), the statistician has access to the vector of ratings $(\phi_k(x))_{1 \leq k \leq n}$.

\item For each user $k$, the statistician has access to $\mu_k = (\phi_k)_*\mu$, the distribution of user $k$'s ratings across all items. Equivalently, the statistician has the cumulative distribution function $F_k$ associated to $\mu_k$ (and hence also the quantile function $F_k^{-1}$).

\end{enumerate}
Recall that for a distribution $\mu$ with cumulative distribution function (CDF) $F_\mu$, the \emph{quantile function} is the generalized inverse function $F_\mu^{-1}(x) := \inf \{ y : F_\mu(y) \geq x\}$. We will define all inverse CDFs in this way so as to ignore invertibility concerns; \cite{embrechts2013note} discusses properties of generalized inverse functions in detail.

The statistician \emph{does not} have access to the functions $\phi_k$, which would otherwise reveal the consensus quality score associated to each item rated by user $k$. In particular, they do not know the original quality score $x$ before transformation into $\phi_k(x)$, only that $\phi_1(x),\dots,\phi_n(x)$ correspond to the same item. Similarly, the statistician lacks access to the distribution $\mu$ of consensus quality scores.

In our analysis of the performance of estimators, we will use $\phi_1,\dots,\phi_n$ to express estimators as functions with domain $([0,1],\mu)$ (the original space of consensus numerical quality scores), even though the functions $\phi_k$ are not actually known to the statistician. For example, the standard average of ratings for an item will be expressed as $A(x) := \frac{1}{n} \sum_{k=1}^n \phi_k(x)$, rather than $A(r_1,\dots,r_n) = \frac{1}{n} \sum_{k=1}^n r_k$ for ratings $r_1,\dots,r_n$ of the unique item with numerical quality score $x$.

Our goal is to provide a function $R : [0,1] \to [0,1]$, using only the data $(\phi_k(x))_{1 \leq k \leq n, x \in [0,1]}$ and $\mu_1,\dots,\mu_n$, such that $R(x)$ recovers the numerical quality score $x$. The loss in estimating $x$ via $R(x)$ will be measured by $|R(x) - x|^2$, averaged over all items; in other words, our estimation error will be $\| R - \id \|_{L^2(\mu)}$, where $\id : [0,1] \to [0,1]$ is the identity function.

So far, our description of the problem has not included any randomness. The only aspect which will be considered ``random'' is how we generate the scale functions $\phi_k$. In our model, we pick users iid from some unknown distribution of users to add to our sample. Keeping in mind the two ``layers'' of distributions (the rating distributions $\mu,\mu_k$ and the overall probability measure $\P$ from which we generate iid user scale functions $\phi_k$), we will always use ``almost everywhere'' or ``a.e.'' when discussing $\mu$-null sets (or $\mu_k$-null sets, etc.) and ``almost surely'' or ``a.s.'' when discussing $\P$-null sets.

\subsection{Identifiability by defining population consensus} \label{identifiability_assumptions}

Thus far, the model is not identifiable, since we have not yet specified any information about the distribution $\mu$ or the personal scale functions $\phi_1,\dots,\phi_n$. The assumptions we make will also serve to separate the problem into two components: the differing \emph{scales} of users' personal rating distributions and the differing \emph{orders} of each user's preferences.

\begin{itemize}

\item[] (A1). (Consensus scale): $\mu$ is the population $W_2$-barycenter of the $\text{Law}(\mu_1)$.

\end{itemize}

Any post-processing of ratings (including doing nothing) induces transformations on users' personal rating distributions. If we want to rescale users' ratings to a common distribution $\nu$ before averaging them (i.e.\ calculating $\frac{1}{n} \sum_{k=1}^n F_\nu^{-1} \circ F_k \circ \phi_k(x)$ for user ratings $\phi_1(x),\dots,\phi_n(x)$ of each item), the canonical choice of $\nu$ would be the measure for which we have to change the user ratings the least; the Wasserstein-2 distance $W_2(\mu_k,\nu)$ quantifies how far we need to move points to rescale user $k$'s ratings to the distribution $\nu$, so the unique distribution $\nu$ which minimizes this cost over all users is the empirical $W_2$-barycenter mean $\wh \mu$. Therefore, the population $W_2$-barycenter can be considered the ``correct common scale'' to view user ratings on, which is precisely the content of (A1).

A variant of assumption (A1) appears functional data analysis literature on phase variation (e.g.\ \cite{10.1214/15-AOS1387}), where random \emph{increasing} homomorphisms $\phi_k : [0,1] \to [0,1]$ satisfy the ``unbiasedness'' condition $\E[\phi_1(x)] = x$ for $\mu$-a.e.\ $x \in [0,1]$. In our setting, users may disagree on item rankings, so the maps $\phi_k$ are not necessarily increasing---this key difference underlies much of Section~\ref{rating_estimation_with_misaligned_preferences} and distinguishes our model from prior work. However, we will see (Proposition~\ref{R0_equals_A_and_is_consistent}) that when users fully agree on preferences (i.e.\ $\phi_k$ is increasing $\P$-a.s.), assumption (A1) reduces to the unbiasedness condition. Thus, our assumption can be viewed as a generalization of unbiased personal scale functions to settings where users disagree on both preferences and scale.

This scale assumption is not sufficient to describe the problem, as this example shows:

\begin{example}
Suppose that $\mu$ is Lebesgue measure on $[0,1]$, and suppose that $\phi_k(x) = 1-x$ $\P$-a.s.\ Then the system satisfies assumption (A1), but there is no way to tell that the items are being rated in the wrong order by all the users. The issue here is that assumption (A1) provides no way to identify \emph{contrarians}, users who report ratings in a different order than the natural order of the consensus numerical quality scores.
\end{example}

Of course, if every user is a contrarian, then no one is. So we need an assumption which allows us to correctly determine the ordering of the quality scores. How should we be able to determine the ordering of the scores of two items with quality scores $x,y \in [0,1]$ (unknown to the statistician) from the users' ratings? Well, if we ask all the users, then the majority opinion should be considered correct. And if we weigh opinions by how strong the users feel about their rating, then a quantity such as $\frac{1}{n} \sum_{k=1}^n \phi_k$ should determine the correct ordering by giving $\frac{1}{n} \sum_{k=1}^n \phi_k(x) < \frac{1}{n} \sum_{k=1}^n \phi_k(y)$ for all sufficiently large $n$ when $x < y$. However, as discussed earlier, the direct average of the ratings is not appropriate because users are rating on different scales. So we should instead hope that $\frac{1}{n} \sum_{k=1}^n F_\mu^{-1} \circ F_k \circ \phi_k(x) < \frac{1}{n} \sum_{k=1}^n F_\mu^{-1} \circ F_k \circ \phi_k(y)$ for all sufficiently large $n$, where $F_k,F_\mu$ are the respective CDFs of $\mu_k,\mu$, respectively. Keeping the strong law of large numbers in mind, which tells us that $\frac{1}{n} \sum_{k=1}^n F_\mu^{-1} \circ F_k \circ \phi_k(x) \to \E[F_\mu^{-1} \circ F_1 \circ \phi_1(x)]$ as $n \to \infty$, we require that $\E[F_\mu^{-1} \circ F_1 \circ \phi_1(x)] < \E[F_\mu^{-1} \circ F_1 \circ \phi_1(y)]$ when $x < y$.

\begin{itemize}

\item[] (A2). (Consensus preference order): The function $g(x) := \E[F_\mu^{-1} \circ F_1 \circ \phi_1(x)]$ is strictly increasing $\mu$-a.e.

\end{itemize}

When $\mu$ is not purely atomic, we replace (A2) with the following technical strengthening.

\begin{itemize}

\item[] (A2'). (Preference order regularity): The functions $g,F_\mu^{-1} \circ F_k \circ \phi_k$ are differentiable with $g' \geq \alpha > 0$ and $(F_\mu^{-1} \circ F_k \circ \phi_k)' \geq -\beta$.

\end{itemize}

\section{The Wasserstein metric and Fr\'echet means} \label{wasserstein_background}

To provide an appropriate estimator for ratings, we would like to first determine a consensus distribution to convert all the user rating distributions $\mu_k$ to.

\begin{example}
One way to average the distributions $\mu_1,\dots,\mu_n$ is to simply use the linear average $\ba \mu := \frac{1}{n} \sum_{k=1}^n \mu_k$. 
%
%
%
A critical issue with this approach is that if users report ratings on a finite set of values, then $\ba \mu$ will report ratings on the same finite set of values. This would lead to many ties in ratings, which is undesirable.
\end{example}

Our choice of average will be based on the Wasserstein 2-distance between rating distributions on $[0,1]$. We now review some basic facts about Wasserstein distances and their associated Fr\'echet means. For a more detailed account of this theory, see e.g.\ \cite{panaretos2020invitation}.

\begin{definition}
The \textbf{Wasserstein 2-distance} (with cost function $c(x,y) = |x-y|^2$) between two probability distributions on $\R$ is 
\begin{align*}
W_2(\mu,\nu) &:= \inf_{\text{couplings $\pi$ of $(\mu,\nu)$}} \sqrt{\int_{\R^2} |x-y|^2 \, d\pi(x,y) } = \inf_{X \sim \mu,Y \sim \nu} \sqrt{\E[|X-Y|^2]}.
\end{align*}
We will use $\mc P_2([0,1])$ for the space of probability measures on $[0,1]$, with this metric.
\end{definition}

For a more general discussion of Wasserstein distances on spaces other than $\R$ with more general cost functions, see \cite{villani2009optimal}.

\begin{lemma}[$W_2$ CDF formula] \label{Wasserstein-CDF-formula}
Let $\mu,\nu \in \mc P_2(\R)$ with respective CDFs $F_\mu, F_\nu$. Then
$$W_2(\mu,\nu) = \sqrt{\int_0^1 |F_\mu^{-1}(x) - F_\nu^{-1}(x)|^2 \, dx}.$$
\end{lemma}

For a proof of Lemma~\ref{Wasserstein-CDF-formula}, see Chapter 1 of \cite{panaretos2020invitation}. The idea is that $F_\nu^{-1} \circ F_\mu$ is the optimal transport map sending $\mu$ to $\nu$. Indeed, the optimal transport map for measures on $\R$ are the increasing functions. Moreover, increasing maps on $\R$ must be of this form.

\begin{lemma} \label{increasing-transport-formula}
Let $\mu,\nu$ be measures on $\R$, and let $f : \R \to \R$ be a function such that $\nu = f_*\mu$. If $f$ is strictly increasing $\mu$-a.e., then $f = F_\nu^{-1} \circ F_\mu$ $\mu$-a.e., where $F_\mu,F_\nu$ are the respective CDFs of $\mu,\nu$.
\end{lemma}

This result is standard. See, for example, Chapter 2 of \cite{santambrogio2015optimal}, which discusses 1-dimensional optimal transport in more detail. We include a short proof in the appendix.

\begin{definition}
Let $(\mc X,\rho)$ be a metric space, and let $x_1,\dots,x_n \in \mc X$. An \textbf{empirical Fr\'echet mean} of $\{x_1,\dots, x_n \}$ is a global minimizer $x \in \mc X$ of the quantity $\frac{1}{n} \sum_{k=1}^n \rho(x_k,x)^2$.
A \textbf{population Fr\'echet mean} of the law of an $\mc X$-valued random variable $X$ is a global minimizer $x \in \mc X$ of the quantity $\E[\rho(X,x)^2]$.
\end{definition}
In general, one must assume that $\E[\rho(X,x_0)^2] < \infty$ for some $x_0$ (and hence for all $x$), but we will only be considering cases where $\rho \leq 1$.

Fr\'echet means in Wasserstein spaces are also sometimes called \emph{barycenters}. In our simple case, the empirical Fr\'echet mean of the personal rating distributions $\mu_k$, using the Wasserstein 2-distance, can be expressed via the inverse CDFs $F_{k}^{-1}$ associated to the $\mu_k$.

\begin{proposition}[$W_2$-barycenter formula] \label{w2-frechet-mean-formula}
Let $\mu_1,\dots,\mu_n$ be measures on $[0,1]$. The empirical $W_2$-barycenter of $\mu_1,\dots,\mu_k$ is the measure $\wh \mu$ with inverse CDF
$$F_{\wh \mu}^{-1} = \frac{1}{n} \sum_{k=1}^n F_k^{-1}.$$
\end{proposition}

We provide a brief proof in the appendix, but a more thorough discussion can be found in \cite{panaretos2020invitation}.

In Section~\ref{rating_estimation_with_misaligned_preferences}, we will make use of Brenier's polar factorization theorem to decouple the contributions of scale and order in ratings estimation. The theorem first appeared in its modern, general form in \cite{brenier1991polar}, but this special case has been long well-known (see e.g.\ \cite{ryff1965orbits}).

\begin{lemma}[Brenier polar factorization] \label{Brenier}
Let $f : [0,1] \to [0,1]$, and let $\mu$ be a probability measure on $[0,1]$. Then there exist a unique nondecreasing function $h : [0,1] \to [0,1]$ and a unique $\mu$-preserving map $\sigma$ such that $f = h \circ \sigma$.
\end{lemma}

Lemma~\ref{increasing-transport-formula} tells us that if $f_*\mu = \nu$, then $h = F_\nu^{-1} \circ F_\mu$. For a detailed discussion of the polar factorization and its generalizations, see Chapter 3 of \cite{villani2021topics}.

\section{Rating estimation with universal preferences} \label{rating_estimation_with_universal_preferences}

In this section, to lay the groundwork for the general case in Section \ref{rating_estimation_with_misaligned_preferences}, we will treat the case where all users agree on all preferences, i.e.\ all users agree whether item $x$ should be rated higher than item $y$, but users still report ratings on individualized scales. We will remove this assumption in favor of (A2) and (A2') in Section~\ref{rating_estimation_with_misaligned_preferences}. Mathematically, the existence of universal preferences is expressed as follows: 
\begin{itemize}

\item[] (A2''). (Universal preferences): $\P$-a.s., $\phi_k$ is increasing $\mu$-a.e.

\end{itemize}
This assumption is unrealistic in practice, and we mainly treat this case because it will elucidate the relationship between our rating estimator and the standard average.

Lemma~\ref{increasing-transport-formula} tells us that with assumption (A2''), $\phi_k$ can be expressed in terms of the CDFs $F_{\mu},F_{k}$ associated to $\mu,\mu_k$, respectively: $\phi_k = F_k^{-1} \circ F_\mu$.

\subsection{The primitive rating estimator} \label{primitive_estimator}

To estimate the quality of items, we can try the following, inspired by \cite{10.1214/15-AOS1387}:

\begin{enumerate}

\item[] Step 1: Calculate the Wasserstein 2-distance Fr\'echet mean $\wh \mu$ of the personal rating distributions $\mu_1,\dots,\mu_n$. Proposition~\ref{w2-frechet-mean-formula} tells us that $F_{\wh \mu}^{-1} = \frac{1}{n} \sum_{k=1}^n F_k^{-1}$.

\item[] Step 2: Estimate the inverse scale functions $\phi_k^{-1}$ by calculating the transport plans $\wh {\phi_k^{-1}}$ from $\mu_k$ to $\wh \mu$. These are increasing, so Lemma~\ref{increasing-transport-formula} shows us that $\wh {\phi_k^{-1}} = F_{\wh \mu}^{-1} \circ F_k$.

\item[] Step 3: Given an item with user ratings $\phi_1(x),\dots,\phi_n(x)$, estimate the quality $x$ via the \emph{primitive rating estimator} $R_0(x) := \frac{1}{n} \sum_{k=1}^n \wh {\phi_k^{-1}} \circ \phi_k(x)$.

\end{enumerate}

Combining these steps, we get the simple formula for the primitive rating estimator:
$$R_0(x) = \frac{1}{n^2} \sum_{k=1}^n \sum_{j=1}^n F_j^{-1} \circ F_k \circ \phi_k(x).$$

\subsection{Asymptotic consistency with universal preferences} \label{asymptotic_consistency_with_universal_preferences}

To discuss consistency of the estimator $R_0$, we first concern ourselves with consistent estimation of $\mu$ via the empirical Fr\'echet mean $\wh \mu$. The consistency of the empirical Fr\'echet mean $\wh \mu$ in estimating the population Fr\'echet mean $\mu$ is well-known (see e.g.\ Chapter 3 of \cite{panaretos2020invitation}) and does not depend on any of the assumptions (A2), (A2'), or (A2'').

\begin{lemma} \label{Wasserstein-barycenter-LLN}
Let $\Lambda$ be a random probability measure in $\mc P_2([0,1])$, and let $\lambda$ be the population $W_2$-barycenter of the law of $\Lambda$. Draw iid samples $\Lambda_1,\Lambda_2,\dots$ according to the law of $\Lambda$, and let $ \lambda_n$ be the empirical $W_2$-barycenter of $\Lambda_1,\dots,\Lambda_n$. Then $W_2(\lambda_n,\lambda) \to 0$ as $n \to \infty$ $\P$-almost surely.
\end{lemma}

Consistency of estimating Wasserstein barycenters was first done in \cite{le2017existence}, which established results in a far more general context than our simple case of $\mc P_2([0,1])$. 


The case of universal preferences is actually much simpler than it may seem at first. Remarkably, it turns out that in this case, the primitive rating estimator is actually the same as $A$, the average of the ratings! Moreover, the loss is naturally expressed in terms of the Wasserstein distance between $\wh \mu$ and $\mu$.

\begin{proposition}[Consistency of $A$, $R_0$ with $L^2$-loss and universal preferences] \label{R0_equals_A_and_is_consistent}
Let $\mu$ be a measure on $[0,1]$, and let $\phi_1,\phi_2,\ldots : [0,1] \to [0,1]$ be iid random functions. Let $\mu_k := (\phi_k)_*\mu$ for each $1\leq k \leq n$, and denote the CDFs of $\mu_k,\mu$ as $F_k,F_\mu$, respectively. Define the assumptions
\begin{itemize}

\item[] (A1). (Consensus preference order): $\mu$ is the population $W_2$-barycenter of $\text{Law}(\mu_1)$.

\item[] (A2''). (Universal preferences): $\P$-a.s., $\phi_k$ is increasing $\mu$-a.e.

\end{itemize}
Define the rating distribution estimator $\wh \mu$ as the Fr\'echet mean of the $\mu_k$ with respect to the Wasserstein 2-distance; i.e.\ $F_{\wh \mu}^{-1} = \frac{1}{n} \sum_{k=1}^n F_k^{-1}$, where $F_{\wh \mu}$ and $F_{k}$ are the CDFs of $\wh \mu, \mu_k$, respectively. Define the average and primitive rating estimator by
$$A := \frac{1}{n} \sum_{k=1}^n \phi_k, \qquad R_0 := \frac{1}{n} \sum_{k=1}^n F_{\wh \mu}^{-1} \circ F_k \circ \phi_k.$$
Then
\begin{enumerate}[label = (\alph*)]

\item Under assumption (A2''), $A = R_0 = F_{\wh \mu}^{-1} \circ F_\mu$ $\mu$-a.e.,

\item $\| F_{\wh \mu}^{-1} \circ F_\mu - \id \|_{L^2(\mu)} = W_2(\wh \mu,\mu)$.
\end{enumerate}
Consequently, under assumptions (A1), (A2''), as $n \to \infty$,
$$\| A - \id \|_{L^2(\mu(m))} = \| R_0 - \id \|_{L^2(\mu(m))} \xrightarrow{\P\text{-a.s.},L^p(\P)} 0 \qquad \text{for all} \quad 1 \leq p < \infty.$$
\end{proposition}

See the appendix for the (very short) proof.

\section{Rating estimation with misaligned preferences} \label{rating_estimation_with_misaligned_preferences}

\subsection{Inconsistency of $A$ and $R_0$ and derivation of $R$}

Moving on to the general case where users can disagree on matters of both preference and scale, the following simple example shows us that $A,R_0$ can differ and that neither of them is a consistent estimator of $\id$ in general. Here, we omit some details from the calculation, but the full details of the calculation may be found in the appendix.

\begin{example} \label{scaling_and_reversing}
The key feature of this example is that the personal scale functions $\phi_k$ will only rescale and reverse preferences. Let $\mu$ be the uniform distribution on $[1/4,3/4]$, and let $\phi_k(x) = \alpha_k(x-1/2) + 1/2$, where the $\alpha_k$ are picked iid from a distribution with $|\alpha_k| \leq 2$ a.s.

We can directly calculate the relevant quantities of interest:
\begin{align*}
W_2(\wh \mu,\mu) = \frac{1}{4 \sqrt 6} \Bigg| \frac{1}{n} \sum_{k=1}^n |\alpha_k| - 1 \Bigg|.
\end{align*}
We know by assumption (A1) and Lemma~\ref{Wasserstein-barycenter-LLN} that $W_2(\wh \mu, \mu) \to 0$, so the strong law of large numbers tells us that $\E[|\alpha_1|] = 1$. Now compare to the loss for the estimators $A$ and $R_0$:
\begin{align*}
\| A - \id \|_{L^2(\mu)} &= \frac{1}{4 \sqrt 6} \Bigg| \frac{1}{n} \sum_{k=1}^n \alpha_k - 1 \Bigg|, \\
\| R_0 - \id \|_{L^2(\mu)} &= \frac{1}{4 \sqrt 6} \Bigg| \Bigg(\frac{1}{n} \sum_{j=1}^n |\alpha_j|\Bigg) \Bigg(\frac{1}{n} \sum_{k=1}^n \sgn(\alpha_k) \Bigg) -1 \Bigg|.
\end{align*}
If $\P(\alpha_1 < 0) > 0$, then neither $\| A - \id \|_{L^2(\mu)}$ nor $\| R_0 - \id \|_{L^2(\mu)}$ converges to 0 as $n \to \infty$. In other words, if users are allowed to disagree at all in their preferences, then neither $A$ nor $R_0$ will accurately recover the numerical quality scores.
\end{example}

The deficiency in using $R_0$ to estimate quality stems from the fact that some users can be \emph{contrarians} with preferences $\phi_k$ which are not increasing, i.e.\ their preference order differs from the consensus preference order. This poses a problem because we try to approximate $\phi_k^{-1}$, which may not be increasing, by increasing functions $\wh {\phi_k^{-1}} := F_{\wh \mu}^{-1} \circ F_k$. This deficiency can be measured by the fact that $R_0$ pushes forward $\mu$ to a measure $\nu$, which may not equal $\wh \mu$ if some of the $\phi_k$ are not increasing.

In other words, because the $\phi_k$ need not always be increasing, $R_0$ can report the scale of the consensus numerical quality scores incorrectly. However, assumption (A2) suggests that $R_0$ asymptotically recovers the ordering of the consensus quality scores, so we can just correct the scale of $R_0$ by transporting the result from $\nu$ to $\wh \mu$. This gives us the \emph{rating estimator}
$$R = F_{\wh \mu}^{-1} \circ F_{\nu} \circ R_0.$$
Recall that $\nu$ is simply the empirical distribution of primitive ratings over all items, so $F_{\nu}$ can be calculated exactly by applying the primitive rating estimator $R_0$ to all items and looking at the resulting distribution of primitive aggregate ratings.

In the appendix, we include a more thorough derivation of the rating estimator, using Brenier's polar factorization theorem. The analysis shows that given desired properties for our estimator, the definition of $R$ is essentially determined.

\begin{remark}
By Proposition~\ref{R0_equals_A_and_is_consistent}, when all users agree on preferences (assumption (A2'')), the rating estimator $R$ equals $A$, the standard average of the ratings. Additionally, whenever all personal rating scales are the same, $R=A$, as well. The benefit of $R$ over $A$ is that it accounts for heterogeneous user rating scales, so when there is no heterogeneity, $A$ is fine.
\end{remark}

\subsection{Asymptotic consistency of the rating estimator} \label{asymptotic_consistency_subsection}

\begin{theorem}[Consistency of $R$ with $L^2$ loss] \label{ratings_are_consistent}



Assuming (A1) and (A2'),
$$\| R - \id \|_{L^2(\mu(m))} \xrightarrow{\P\text{-a.s.},L^p(\P)} 0 \qquad \text{for all} \quad 1 \leq p < \infty$$
as $n \to \infty$. If $\mu$ is purely atomic, then assumption (A2') can be weakened to (A2).



\end{theorem}

\begin{remark}
When there are only finitely many items to rate, the measure $\mu$ will be purely atomic, in which case we only need the minimal assumptions (A1), (A2).
\end{remark}

We leave the proof to the appendix, which involves proofs of the following two Glivenko--Cantelli-type theorems.

\begin{theorem}[Glivenko--Cantelli-type theorem for quantile functions] \label{inverse_CDF_GC}
Let $\Lambda$ be a random measure on $[0,1]$, and let $\lambda$ be the population $W_2$-barycenter of the law of $\Lambda$. Draw iid samples $\Lambda_1,\Lambda_2,\dots$ according to the law of $\Lambda$, and let $ \lambda_n$ be the empirical $W_2$-barycenter of $\Lambda_1,\dots,\Lambda_n$. Denote the CDFs of $ \lambda_n$ and $\lambda$ as $F_{ \lambda_n}, F_\lambda$, respectively. 
Then
$$\sup_{x \in [0,1]} | F_{ \lambda_n}^{-1}(x) - F_{\lambda}^{-1}(x)| \to 0$$
almost surely as $n \to \infty$. Moreover, if $\supp \lambda = [0,1]$, then
$$\E \squa{\sup_{x \in [0,1]} | F_{ \lambda_n}^{-1}(x) - F_{\lambda}^{-1}(x)|} \leq 4 \sqrt{\frac{\log n}{n}} \qquad \forall n \geq 3.$$
\end{theorem}

\begin{theorem}[Glivenko--Cantelli-type theorem for CDFs] \label{wass-gc}
With the same notation as Theorem~\ref{inverse_CDF_GC}, assume that $\lambda$ has no atoms and that $\supp \Lambda = [0,1]$ a.s. Then
$$\sup_{x \in [0,1]} | F_{ \lambda_n}(x) - F_{\lambda}(x)| \to 0$$
almost surely as $n \to \infty$.
\end{theorem}

\begin{remark}[Order-scale decomposition for $L^2$ loss]
The technique in our proof can be used to analyze any estimator $E$ of consensus ratings. Brenier's polar factorization characterizes the $\mu$-preserving map $\sigma$ (associated to $E$) as the $L^2$-projection of $E$ onto the set of $\mu$-preserving maps. In particular,
\begin{align*}
\| E - \id \|_{L^2(\mu)}^2 &= \| E - \sigma \|_{L^2(\mu)}^2 + \| \sigma - \id \|_{L^2(\mu)}^2 \\
&= W_2(E_*\mu,\mu)^2 + \| \sigma - \id \|_{L^2(\mu)}^2.
\end{align*}
\end{remark}

\begin{theorem}[Rates of convergence for the rating estimator] \label{rating_rates_of_convergence}
With the same notation as Theorem~\ref{ratings_are_consistent},
\begin{enumerate}[label = (\alph*)]

\item If $\mu$ is purely atomic with $M$ atoms, assuming (A1) and (A2), then for all $n \geq 1$,
\begin{align*}
\E[\| R - \id \|_{L^2(\mu)}] &\leq \frac{1}{2\sqrt n} + 4M e^{-n\delta^2/8} = O \paren{\frac{1}{\sqrt n}},
\end{align*}
where $\delta := \min \{ g(y) - g(x) : x,y \text{ are atoms of $\mu$}, x < y  \}$.

\item If $\supp \mu = [0,1]$, assuming (A1) and (A2'), then for all $n \geq 3$,
\begin{align*}
\E[\| R - \id \|_{L^2(\mu)}] &\leq \paren{\frac{3C}{\alpha} + \frac{11}{2}} \sqrt{\frac{\log n}{n}} = O \paren{\sqrt{\frac{\log n}{n}}},
\end{align*}
where $C$ is a constant depending only on $\beta$.
\end{enumerate}
\end{theorem}


\begin{remark}
The rate of convergence $1/\sqrt n$ is optimal. Recall that when all users agree on preference order, Proposition~\ref{R0_equals_A_and_is_consistent} tells us that $R = A$, so we are estimating $x = \E[\phi_1(x)]$ by $\frac{1}{n} \sum_{k=1}^n \phi_k(x)$, which follows the central limit theorem rate of $1/\sqrt n$. 
\end{remark}

\section{Rating with incomplete user data} \label{incomplete-data}

In many settings, not every user has rated every item. We now show that the rating estimator still works, given assumptions regarding how the data is missing at random.



Here is some notation. Denote $N_x$ as the (possibly random) set of users rating the item with quality score $x$, and let $N = \{ N_x : x \in \supp \mu \}$ be the collective information of ``who rates what.'' Rather than $\mu_1,\dots,\mu_n$, the statistician has access to $N$, along with the \emph{empirical personal rating distributions} $\wt \mu_1,\dots,\wt \mu_n$, represented by the empirical CDFs $\wt F_1,\dots,\wt F_n$. The empirical Wasserstein-2 Fr\'echet mean $\wt \mu$ is given by $F_{\wt \mu}^{-1} := \frac{1}{n} \sum_{k=1}^n \wt F_k^{-1}$, and our adjusted rating estimator is
$$\wt R := F_{\wt \mu}^{-1} \circ F_{(\wt R_0)_*\mu} \circ \wt R_0, \qquad \text{where} \quad  \wt R_0(x) := \frac{1}{|N_x|} \sum_{k \in N_x} F_{\wt \mu}^{-1} \circ \wt F_k \circ \phi_k(x).$$
It will be mathematically useful to maintain our notion of the personal scale functions $\phi_k$ and to assume that they are defined $\mu$-a.e., regardless of whether or not the corresponding user has rated each item. We interpret $\phi_k(x)$ as the \emph{hypothetical rating} of user $k$ for the item with quality score $x$. Two new challenges arise in our analysis:

\begin{enumerate}[label = (\alph*)]

\item Do the items user $k$ rates give an accurate depiction of their personal rating scale $\mu_k$?

\item Have enough users rated each item for us to accurately understand its quality?

\end{enumerate}

To address the second concern, we will necessarily need to work in the regime where the number of ratings for each item goes to infinity. To address the first concern, we will need assumptions regarding which items each user rates. In reality, there are a number of different ways in which ``who rates what'' is chosen, and one set of assumptions may not apply to all contexts. Here, we aim to provide relatively general assumptions to serve as a proof of concept in justifying the robustness of our rating estimator.

\begin{theorem} \label{incomplete-consistency-thm}
Let $\mu$ be purely atomic with $M$ atoms, and assume that

\begin{itemize}

\item[] (A1). (Consensus scale): $\mu$ is the population $W_2$-Fr\'echet mean of $\text{Law}(\wt \mu_1)$.

\item[] ({{A2'}'}'). (Consensus preference order with incomplete data): The function $h(x) := \E[F_\mu^{-1} \circ \wt F_1 \circ \phi_1(x)]$ is increasing.

\item[] (A3). (Balanced choices): $\E[\wt F_1^{-1}(x) \mid \phi_1] = F_1^{-1}(x)$ for all $x \in [0,1]$.

\item[] (A4). (Item choices don't influence profiles): The collection $N := \{ N_x : x \in \supp \mu\}$ of ``who rates what'' is independent of $\{ \phi_k : k \geq 1\}$.

\end{itemize}
Then
$$\E[\| \wt R - \id \|_{L^2(\mu)}] \leq \frac{1}{2\sqrt n} + 6.5 \sum_{\text{atoms $x$}} \E[e^{-|N_x|\delta^2/128}],$$
where $\delta := \min \{ h(y) - h(x) : \text{$x,y$ are atoms of $\mu$, $x < y$}\}$. Consequently, if for all atoms $x$ of $\mu$, $|N_x| \geq \frac{64}{\delta^2} \log n$, then
\begin{align*}
\E[\| \wt R - \id \|_{L^2(\mu)}] &\leq \frac{7M}{\sqrt n} = O\paren{\frac{1}{\sqrt n}}.
\end{align*}
\end{theorem}

The constants are not optimized and have been weakened to simplify the presentation. 

\begin{remark}
Together, assumptions (A1) and (A3) imply that $\mu$ is the $W_2$-Fr\'echet mean of the law of $\wt \mu_1$, as $\E[\wt F_1^{-1}(x)] = F_\mu^{-1}(x)$ for all $x \in [0,1]$. One can prove (with almost the exact same proof, verbatim) the same theorem below by replacing (A1) and (A3) by the assumption that $\mu$ is the $W_2$-Fr\'echet mean of the law of $\wt \mu_1$. However, we find that the assumptions (A1), (A2), (A3) more closely match how one conceptualizes rating systems.
\end{remark}

\begin{remark}
In some applications, (A3) may not be exactly satisfied. A careful modification of our proof will show that (A3) can be replaced with a suitable asymptotic statement as the number of items each user rates goes to infinity. However, we use (A3) to avoid overcomplicating the proof.
\end{remark}

\section{Rankings and Kendall's W} \label{ranking-section}

A common alternative type of data to numerical ratings is rankings, i.e.\ users provide a preference list for items. This can be viewed as a special case of ratings by viewing a rank of $k$ out of $M$ items as the numerical rating $1- (k-1)/M$, which is its quantile among the items. In this setting, the distances between users' successive numerical ratings are uniformly set to $1/M$, and users all have the same personal rating scale, $\mu_1 = \mu_2 = \cdots = \mu_n = \mu = \text{Uniform}(\{1/M,2/M,\dots,1 \})$. Consequently, the primitive rating estimator $R_0$ equals the average, $A$. So all three of $A$, $R_0$, and $R$ give the same consensus preference order, with the rating estimator $R$ outputting the order in the same equally spaced format as the input, i.e.\ with distribution $\mu$. The main purpose of this section is to point out this relationship between ratings and rankings and to show by way of example that methods of studying rankings can be generalized to study ratings, as well.

The discrepancy between $A$ and $R$ in this context has been studied as Kendall's W statistic \cite{kendall1939problem}. Kendall's W is a widely used (e.g.\ \cite{field2005kendall, legendre2005species, gearhart2013use}) statistic for measuring the degree of agreement on preference order when multiple users rank items. Given a matrix $(m_{i,j})_{1 \leq i \leq n, 1 \leq j \leq M}$ of rank $i$ of item $j$, let $\ba m_j := \frac{1}{n} \sum_{i=1}^n m_{i,j}$ denote the average rank given to item $j$, and let $\ba m := \frac{1}{M} \sum_{i=1}^n \ba m_i$ be the grand mean of the matrix. Kendall's W is given by
\begin{align*}
W &:= \frac{\frac{1}{M}\sum_{j=1}^M (\ba m_j - \ba m)^2}{ \frac{1}{12}(M^2-1)} = \frac{\Var(A_*\mu)}{\Var(\mu)} = \frac{\Var(A_*\mu)}{\frac{1}{n} \sum_{i=1}^n \Var(\mu_i)}.
\end{align*}
(Here, we mean the variance associated to these measures, rather than with respect to the randomness generated by the underlying probability measure $\P$.)

Intuitively, if there is disagreement among users, the differing opinions will ``destructively interfere,'' leading the average ranks to tend away from the extreme values $1/M, 1$. Accordingly, $0 \leq W \leq 1$, with the equality cases being when all items are tied after tallying the rankings (no concordance in preferences) and when all users submit the exact same rankings for all items. Following this intuition, we can propose analogous statistics for measuring the concordance in user rating scales and in user ratings overall:
$$W_{\text{scale}} := \frac{\Var(R_*\mu)}{\frac{1}{n} \sum_{i=1}^n \Var(\mu_i)}, \qquad W_{\text{ratings}} := \frac{\Var((R_0)_*\mu)}{\frac{1}{n} \sum_{i=1}^n \Var(\mu_i)}.$$
In general, $R_*\mu$ is $\wh \mu$, the empirical $W_2$-barycenter of $\mu_1,\dots,\mu_n$. Since $R_0 = A$ for rankings, $W_{\on{ratings}}$ generalizes Kendall's $W$ from the special case of rankings to general ratings.

We record formulas for these quantities, although in practice, it may be faster to calculate $W_{\text{scale}},W_{\text{ratings}}$ from $R_*\mu, (R_0)_*\mu$ after calculating $R,R_0$, respectively for all items.

\begin{proposition}[Concordance statistic formulas] \label{concordance-formulas}
Given a matrix $(r_{i,j})_{1 \leq i \leq n, 1 \leq j \leq M}$ of rating $i$ of item $j$, let $\ba r_i := \frac{1}{M} \sum_{j=1}^M r_j$ be user $i$'s rating mean and $\sigma_i^2 := \frac{1}{M} \sum_{j=1}^M (r_{i,j} - \ba r_i)^2$ user $i$'s rating variance.
\begin{enumerate}[label = (\alph*)]

\item
\begin{align*}
W_{\on{scale}} &= \frac{\sum_{i=1}^n \sum_{k=1}^n \Cov(X_i,X_k)}{n \sum_{i=1}^n \Var(X_i)} \\
&= \frac{\sum_{i=1}^n \sum_{k=1}^n (\sum_{j=1}^M r_{i,(j)}r_{k,(j)}) - M\ba r_i \ba r_k }{n M \sum_{i=1}^n \sigma_i^2}
\end{align*}
where $X_i \sim \mu_i$ and the $X_i$ are monotone coupled (i.e.\ $X_k = F_k^{-1} \circ F_i(X_i)$ for each $i,k$). Here, $r_{i,(j)}$ means the $j$-th lowest rating from user $i$.

\item
\begin{align*}
W_{\on{ratings}} &= \frac{\sum_{i=1}^n \sum_{k=1}^n \sum_{\ell=1}^n \sum_{p=1}^n \Cov(Y_{i,k},Y_{\ell,p}) }{n^3 \sum_{i=1}^n \Var(X_i)} \\
&= \frac{\sum_{i=1}^n \sum_{k=1}^n \sum_{\ell=1}^n \sum_{p=1}^n (\sum_{j=1}^M r_{k,(\on{ind}(r_{i,j}))} r_{p,(\on{ind}(r_{\ell,j}))}) - M\ba r_i \ba r_k }{n^3 M \sum_{i=1}^n \sigma_i^2},
\end{align*}
where $Y_{i,k} := F_k^{-1} \circ F_i(r_{i,J})$ and $J \sim \on{Uniform}(\{1,\dots,M\})$. Here, $\on{ind}(r_{i,j})$ denotes the index of $r_{i,j}$ among user $i$'s sorted ratings $r_{i,(1)},\dots,r_{i,(M)}$, and $r_{k,(\ind(r_{i,j}))}$ is the $\ind(r_{i,j})$-th smallest rating given by user $k$.

\end{enumerate}
\end{proposition}


As with Kendall's $W$, these statistics give values between 0 and 1, with values closer to 1 expressing higher concordance in user rating scales and user rating profiles, respectively.

\begin{proposition} \label{concordance-bounds}
Assume that $\mu_1,\dots,\mu_n$ are not all point masses (so the denominators of $W_{\text{scale}}, W_{\text{ratings}}$ are nonzero).
\begin{enumerate}[label = (\alph*)]

\item $0 \leq W_{\text{scale}} \leq 1.$ Equality on the left never occurs. Equality on the right occurs iff all $\mu_i$ are the same, up to translation by a constant.

\item $0 \leq W_{\text{ratings}} \leq 1.$ Equality on the left occurs iff user opinions on every item are equally balanced, i.e.\ $R_0$ is constant. Equality on the right occurs iff all user ratings are the same, up to translation by a constant, i.e.\ $\phi_i = \phi_\ell + a_{i,\ell}$ for all $1 \leq i, \ell \leq M$.

\end{enumerate}
\end{proposition}


To provide a sense of calibration for these statistics, we calculate the values for Example~\ref{scaling_and_reversing}. We record only the results here, but the full calculations can be found in the appendix.

\begin{example}
Let $\mu$ be any distribution on $[0,1]$ which is symmetric about $1/2$ (besides the point mass at $1/2$), and let $\phi_k(x) := \alpha_k(x-1/2) + 1/2$ for some iid $\alpha_k$ with $\E[|\alpha_1|] = 1$ and $\P(\alpha_1 > 0) > 1/2$. Then
$$W_{\on{scale}} = \frac{(\frac{1}{n} \sum_{k=1}^n |\alpha_k|)^2}{\frac{1}{n} \sum_{k=1}^n \alpha_k^2} \xrightarrow{a.s., n \to \infty} \frac{1}{\E[\alpha_1^2]},$$
$$W_{\on{ratings}} = \frac{(\frac{1}{n} \sum_{k=1}^n |\alpha_k|)^2 (\frac{1}{n} \sum_{k=1}^n \sgn(\alpha_k))^2}{\frac{1}{n} \sum_{k=1}^n \alpha_k^2} \xrightarrow{a.s., n \to \infty} \frac{(\E[\sgn(\alpha_1)])^2}{\E[\alpha_1^2]}.$$
Since the distribution of $\alpha_k$ is subject to the constraint that $\E[|\alpha_k|] = 1$, $W_{\on{scale}}$ decreases as $\Var(|\alpha_1|)$ increases. Similarly, $W_{\on{ratings}}$ decreases as $\Var(|\alpha_1|)$ increases, with a multiplicative factor which decreases the value as the population agrees less on preference, i.e. when $\P(\alpha_1 < 0)$ grows closer to $1/2$.
\end{example}

\section{Numerical experiments with real-world rating data} \label{application_to_real_data}

We have applied the estimators analyzed in the previous sections to real and simulated data, including datasets from the Netflix prize challenge, FoodPanda Bangladesh, MyAnimeList, and Yelp.

\subsection{Challenges with benchmarking rating aggregation methods}

Rating data are subjective, so there is no single “correct’’ notion of aggregate ratings and no verifiable “ground truth’’ item-quality scores. As a result, we cannot evaluate aggregation methods using relevance-based measures such as normalized discounted cumulative gain \cite{jarvelin2002cumulated} or mean average precision \cite{harman1992evaluation}. Ordinal comparison methods such as Kendall’s $\tau$ \cite{665905b2-6123-3642-832e-05dbc1f48979}, which count the number of matching pairwise preferences, are also misaligned with rating aggregation, whose goal is to capture not just pairwise preferences but \emph{the degree to which users prefer one item to another}. Accordingly, to enable a \emph{quantitative} comparison of aggregation methods, we evaluate their performance on a downstream task---such as rating prediction (even though aggregate ratings are not explicitly intended to predict user ratings)---as a proxy for alignment with the population’s preferences.

\subsection{Overall distributions and performance on prediction tasks}

The basic prediction task we employ is motivated by the idea that if a global ranking of items aligns with the population’s preferences, then, for example, an item at the 70th percentile of the ranking should on average receive ratings near the 70th percentile of users’ rating distributions. Formally, given a ranking $\mathbf{r} = (r_1,\dots,r_M)$ of the $M$ items, we predict that user $i$ would rate item $j$ as $F_i^{-1}(\frac{M + 1 - r_j}{M})$.

We assess performance via 5-fold cross-validation. In each fold, we learn an aggregate rating (and its induced ranking) from 80\% of the users. For the remaining 20\% of test users, we randomly hide 20\% of their ratings (using the same hidden entries across methods), estimate each user’s personal rating distribution from the remaining 80\% of their ratings, and then predict the hidden ratings. We finally average the prediction error over these hidden ratings. We include the average result of 5 random rankings for comparison. The code to reproduce these calculations (or to apply these estimators to your own rating data) is available at \url{https://github.com/pillowmath/Rating-Estimator}.

\begin{figure}
\begin{center}
\includegraphics[scale = 0.33]{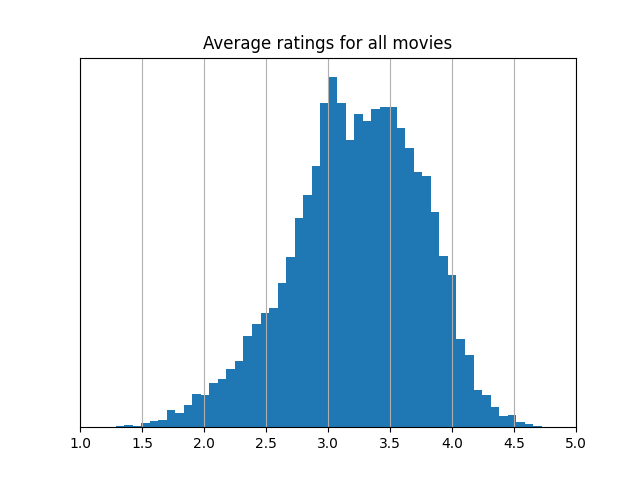} \includegraphics[scale = 0.33]{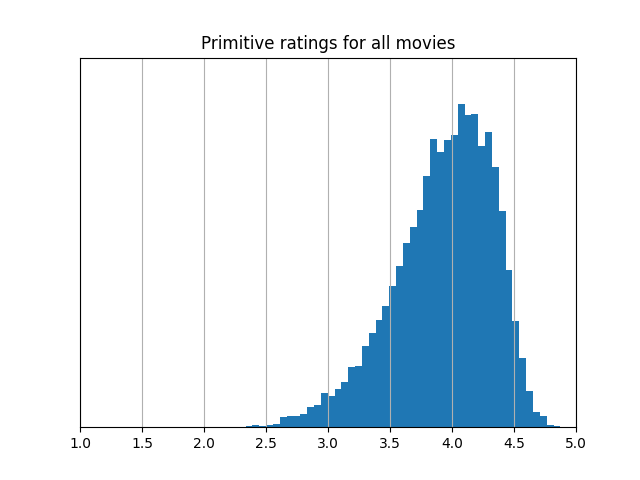} \includegraphics[scale = 0.33]{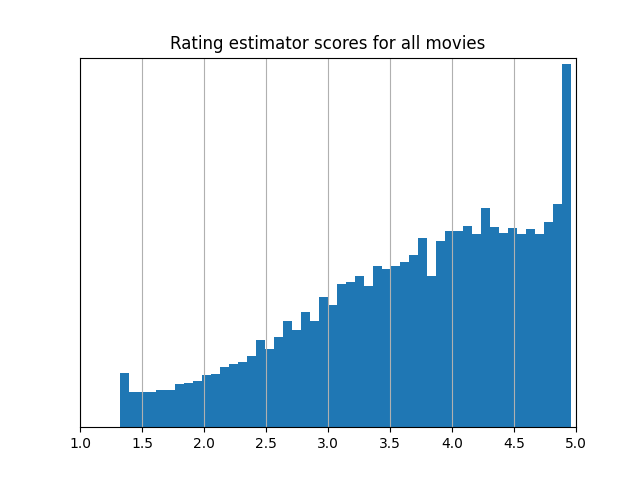}
\end{center}
\caption{The distributions $A_*\mu$ (left), $(R_0)_*\mu$ (middle), and $\wh \mu$ (right) for the Netflix prize dataset, generated by applying the average $A$, the primitive rating estimator $R_0$, and the rating estimator $R$, respectively.}
\label{netflix-distributions}
\end{figure}

\textbf{Netflix prize challenge.} The Netflix prize dataset \cite{bennett2007netflix} is a well-known dataset with over 480,000 users rating over 17,000 movies on a 1 to 5 star integer scale. Figure~\ref{netflix-distributions} shows the resulting rating distributions (over all rated movies in the sample) when using the estimators $A$, $R_0$, and $R$, respectively. Since not every user in the sample has rated every movie, we modified the rating estimators as discussed in Section~\ref{incomplete-data}. Figure~\ref{netflix-distributions} shows the resulting rating distributions (over all rated movies in the sample) when using the estimators $A$, $R_0$, and $R$, respectively, and Figure~\ref{netflix-predictions} shows the performance on the prediction task. The rating estimator scores are less concentrated, and the spikes at the ends of the distribution indicate that users give many ratings at the top or the bottom of their lists. To measure the level of agreement among personal rating scales and among all user rating profiles, we calculate our W statistics, as proposed in Section~\ref{ranking-section}, for this dataset:
$$W_{\text{scale}} = 0.789, \qquad W_{\text{ratings}} = 0.149.$$
From the relatively high value of $W_{\text{scale}}$ and relatively low value of $W_{\text{ratings}}$, we surmise that Netflix users' personal rating scales tend to be relatively similar in scale, but user ratings differ heavily in their preference orders for items.

\begin{figure}
\begin{center}
$$ \begin{array}{|c|c|c|c|}
\hline
\textbf{Netflix Ratings} & \text{MAE (SE)} & \text{MSE (SE)} & \text{Exact Match Rate (SE)} \\
\hline
\text{Average} & 0.831 \: (0.0002) & 1.445 \: (0.0006) & 0.407 \: (0.0001) \\
\hline
\text{Rating Estimator} & 0.778 \: (0.0002) & 1.293 \: (0.0005) & 0.426 \: (0.0001) \\
\hline
\text{Avg of 5 Random Rankings} & 1.037 \: (0.0001) & 1.982 \: (0.0003) & 0.322 \: (0.00005) \\
\hline
\end{array}$$
\end{center}
\caption{Performance of estimators on a simple prediction task for the Netflix dataset. The rating estimator consistently outperforms the standard average.}
\label{netflix-predictions}
\end{figure}

\begin{figure}
\begin{center}
\includegraphics[scale = 0.5]{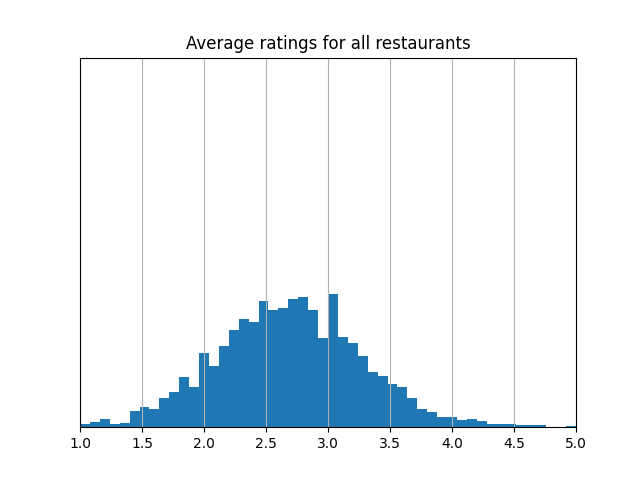} \includegraphics[scale = 0.5]{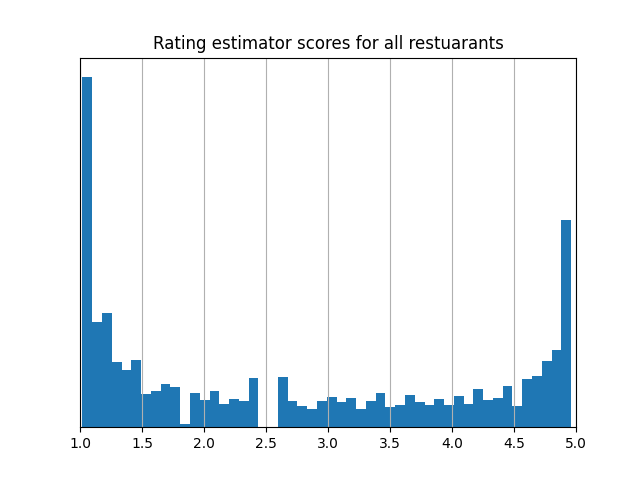}
\end{center}
\caption{The distributions $A_*\mu$ (left), $\wh \mu$ (right) for the Foodpanda restaurant dataset, generated by applying the average $A$ and the rating estimator $R$, respectively.}
\label{Foodpanda-distributions}
\end{figure}

\textbf{Foodpanda restaurant ratings.} Foodpanda is an online food delivery service that allows users to rate and review restaurants. This dataset \cite{deb2024bangladesh} of ratings from Foodpanda's Bangladeshi users has over 300,000 1 to 5 star ratings from over 3,800 users of over 3,600 restaurants. To improve the data quality, we removed users and restaurants with fewer than 10 ratings associated to them. Figure~\ref{Foodpanda-distributions} shows the resulting rating distributions, and Figure~\ref{Foodpanda-predictions} shows the performance on the prediction task.

Average ratings form a bell curve, but the rating estimator is strongly bimodal, with peaks near 1 and 5 stars. Low estimator scores typically correspond to restaurants near the bottom of users’ lists, and high scores to those near the top, making it easier to separate good from bad options. In contrast, average ratings are compressed into a narrow band (around 2–3.5), so two restaurants with similar averages may still differ a lot in quality. Because the estimator reports aggregate ratings on a similar scale to users’ ratings, the results are more interpretable than a simple average; for instance, users are likely to give a 5 star aggregate rated restaurant a 5 star rating (and similarly for 1 star rated restaurants).

The $W$ statistics for this dataset are
$$W_{\text{scale}} = 0.818, \qquad W_{\text{ratings}} = 0.077,$$
indicating that there may be even less preference agreement on restaurants than on movies in the Netflix dataset.

\begin{figure}
\begin{center}
$$ \begin{array}{|c|c|c|c|}
\hline
\textbf{Foodpanda Ratings} & \text{MAE (SE)} & \text{MSE (SE)} & \text{Exact Match Rate (SE)} \\
\hline
\text{Average} & 1.489 \: (0.007) & 4.104 \: (0.028) & 0.3245 \: (0.0024) \\
\hline
\text{Rating Estimator} & 1.495 \: (0.007) & 4.128 \: (0.028) & 0.3222 \: (0.0024) \\
\hline
\text{Avg of 5 Random Rankings} & 1.794 \: (0.003) & 5.435 \: (0.014) & 0.2769 \: (0.0010) \\
\hline
\end{array}$$
\end{center}
\caption{Performance of estimators on a simple prediction task for the Foodpanda dataset. The rating estimator performs about as well as the standard average.}
\label{Foodpanda-predictions}
\end{figure}

\begin{figure}
\begin{center}
\includegraphics[scale = 0.33]{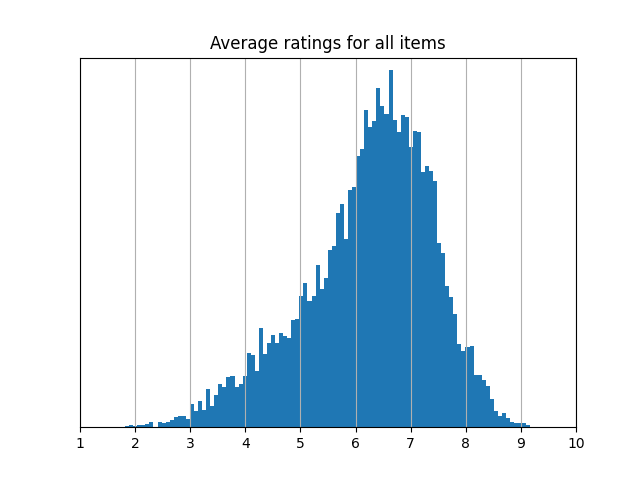} \includegraphics[scale = 0.33]{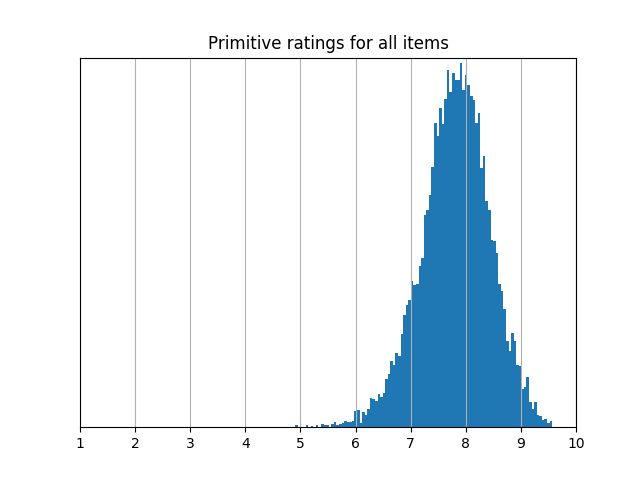} \includegraphics[scale = 0.33]{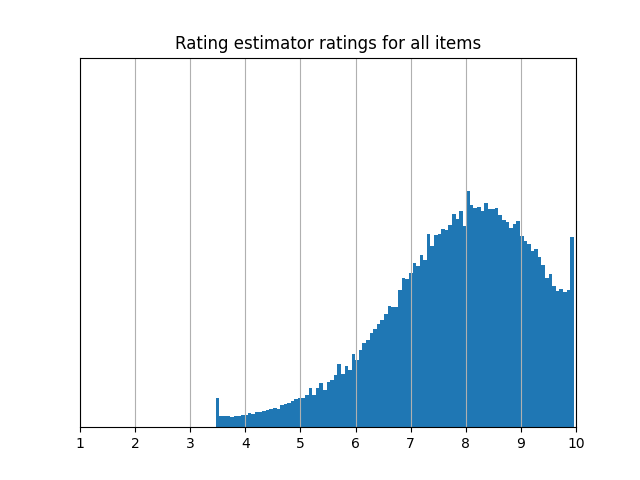}
\end{center}
\caption{The distributions $A_*\mu$ (left), $(R_0)_*\mu$ (middle), and $\wh \mu$ (right) for the MyAnimeList dataset, generated by applying the average $A$, the primitive rating estimator $R_0$, and the rating estimator $R$, respectively.}
\label{mal-distributions}
\end{figure}

\textbf{MyAnimeList ratings.} \url{https://www.MyAnimeList.net} is a rating website where users rate Japanese animation (anime). The dataset \cite{mal2018dataset} contains over 46 million ratings from over 300,000 users on over 14,000 anime, on a range of 1 to 10. Figure~\ref{mal-distributions} shows the resulting rating distributions, and Figure~\ref{mal-predictions} shows the performance on the prediction task. The rating estimator scores are less concentrated, and the spikes at the ends of the distribution indicate that users give many ratings at the top or the bottom of their lists. Our W statistics for this dataset are
$$W_{\text{scale}} = 0.814, \qquad W_{\text{ratings}} = 0.178.$$
The heterogeneity of user rating scales and preference orders seem to be comparable to those of the Netflix dataset.

\begin{figure}
\begin{center}
$$ \begin{array}{|c|c|c|c|}
\hline
\textbf{MyAnimeList Ratings} & \text{MAE (SE)} & \text{MSE (SE)} & \text{Exact Match Rate (SE)} \\
\hline
\text{Average} & 1.552 \: (0.0004) & 4.203 \: (0.002) & 0.219 \: (0.0001) \\
\hline
\text{Rating Estimator} & 1.353 \: (0.0004) & 3.407 \: (0.002) & 0.263 \: (0.0001) \\
\hline
\text{Avg of 5 Random Rankings} & 1.542 \: (0.0002) & 4.416 \: (0.001) & 0.238 \: (0.0001) \\
\hline
\end{array}$$
\end{center}
\caption{Performance of estimators on a simple prediction task for the MyAnimeList dataset. The rating estimator significantly outperforms the standard average.}
\label{mal-predictions}
\end{figure}

\begin{figure}
\begin{center}
\includegraphics[scale = 0.33]{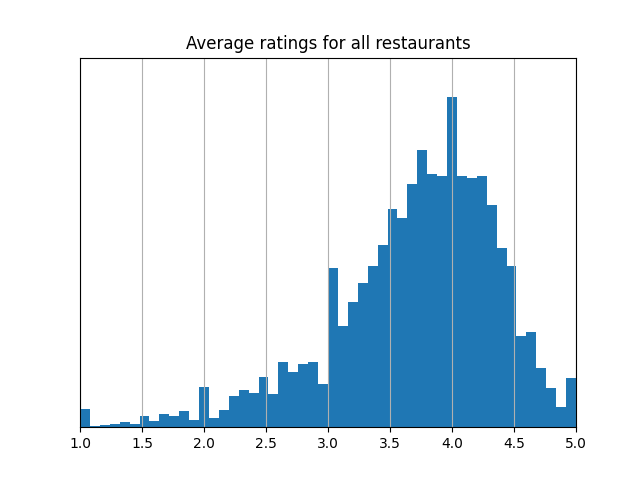} \includegraphics[scale = 0.33]{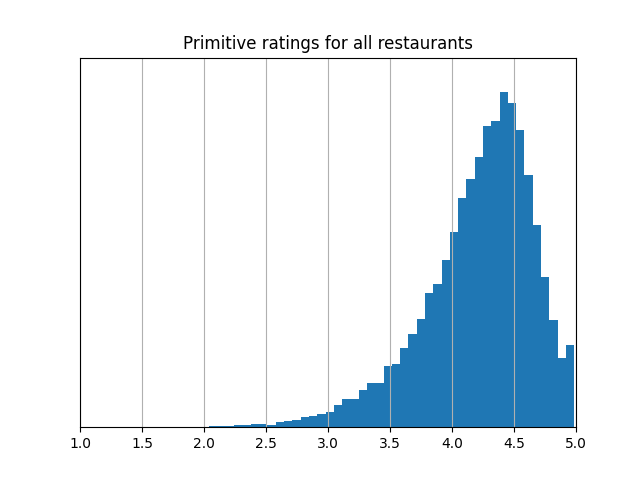} \includegraphics[scale = 0.33]{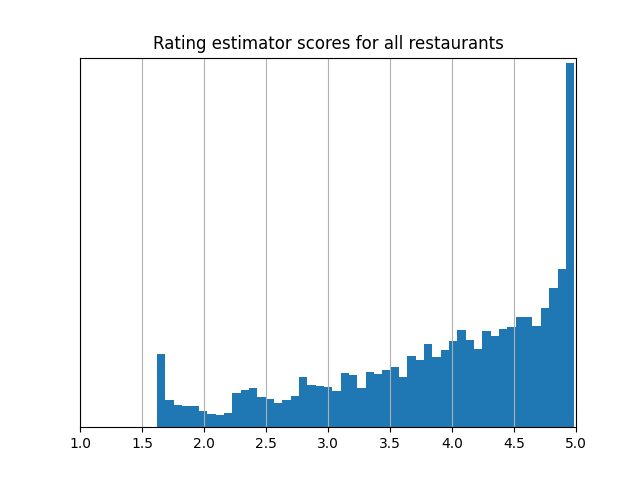}
\end{center}
\caption{The distributions $A_*\mu$ (left), $(R_0)_*\mu$ (middle), and $\wh \mu$ (right) for the Yelp dataset, generated by applying the average $A$, the primitive rating estimator $R_0$, and the rating estimator $R$, respectively.}
\label{yelp-distributions}
\end{figure}

\textbf{Yelp restaurant ratings.} We have also applied the rating estimator to the Yelp open dataset \cite{yelp_kaggle_dataset}, which contains 1 to 5 star ratings of over 60,000 restaurants from over 2 million users. To improve the data quality, we removed users and restaurants with fewer than 15 ratings (an arbitrary cutoff) associated to them. Figure~\ref{yelp-distributions} shows the resulting rating distributions, and Figure~\ref{yelp-predictions} shows the performance on the prediction task. Our W statistics for this dataset are
$$W_{\text{scale}} = 0.737, \qquad W_{\text{ratings}} = 0.172,$$
indicating less scale agreement but more preference agreement than in the Foodpanda restaurant ratings.

\begin{figure}
\begin{center}
$$ \begin{array}{|c|c|c|c|}
\hline
\textbf{Yelp Ratings} & \text{MAE (SE)} & \text{MSE (SE)} & \text{Exact Match Rate (SE)} \\
\hline
\text{Average} & 0.7541 \: (0.0020) & 1.3360 \: (0.0058) & 0.4506 \: (0.0011) \\
\hline
\text{Rating Estimator} & 0.7552 \: (0.0020) & 1.3316 \: (0.0058) & 0.4477 \: (0.0011) \\
\hline
\text{Avg of 5 Random Rankings} & 1.0349 \: (0.0010) & 2.1315 \: (0.0035) & 0.3424 \: (0.0005) \\
\hline
\end{array}$$
\end{center}
\caption{Performance of estimators on a simple prediction task for the Yelp dataset. The rating estimator performs about as well as the standard average.}
\label{yelp-predictions}
\end{figure}

\subsection{Application on simulated data}

One must also take care when simulating ratings, as the choice of simulation model can bias the results towards one aggregation method. Moreover, ratings derived from rankings are not sufficient for comparison, as there is no heterogeneity in personal rating scales (and therefore $R = A$); even with iid Bernoulli($p$) omitted ratings, one can check that $F^{-1}_{\widehat \mu}$ concentrates around the discrete uniform quantile function, and we still have $R \approx A$. We analyze two simulations, based on Example~\ref{scaling_and_reversing}'s ``scale and reverse'' model and a simple ``true quality plus noise'' model, on which the use of the standard average is generally predicated.

\begin{figure}
\begin{center}
\includegraphics[scale = 0.47]{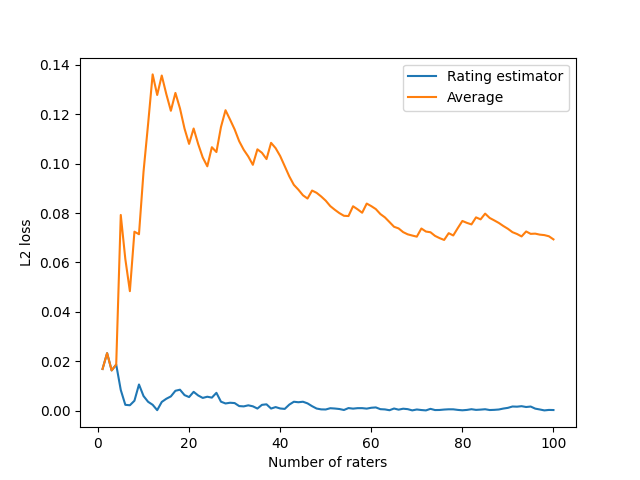} \includegraphics[scale = 0.4275]{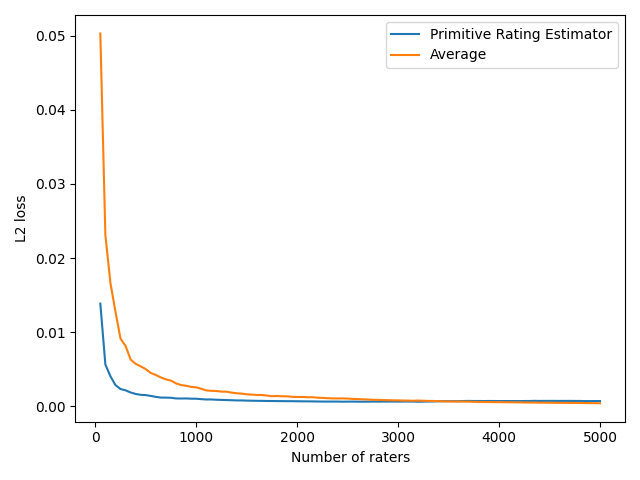}
\end{center}
\caption{$L^2$ losses for our estimators vs the average on simulated data from a scale and reverse model (left) and a quality plus noise model (right).}
\label{simulation-plots}
\end{figure}

In Example~\ref{scaling_and_reversing}, $\mu$ is uniform on $[1/4,3/4]$, and $\phi_k = \varepsilon_k Z_k (x-1/2) + 1/2$, where $Z_k \sim N(1,1/16)$, $\varepsilon_k$ equals $+1$ with probability $3/4$ and $-1$ with probability $1/4$, and all the $Z_k,\varepsilon_k$ are mutually independent. In other words, we apply a Gaussian random scaling and a biased random reversal of preferences. Figure~\ref{simulation-plots}(a) plots the $L^2$ loss.
As predicted by Example~\ref{scaling_and_reversing} and Theorem~\ref{ratings_are_consistent}, $R$ significantly outperforms the average for the $L^2$ loss.

The quality plus noise model generates ratings from user $1 \leq i \leq n$ of item $1 \leq j \leq M$ as $r_{i,j} \coloneqq \mu_i + \sigma_i (q_j + \varepsilon_{i,j})$, where $\E[\varepsilon_{i,j}] = 0$ and $q_j$ are the quality scores. Without loss of generality, we may assume $\E[\mu_i] = 0$ and $\E[\sigma_i] = 1$. If we also assume $\varepsilon_{i,j},\sigma_i,\mu_i$ are independent, then $\ba r_j \coloneqq \frac{1}{n} \sum_{i=1}^n r_{i,j}$ in an unbiased estimator of $q_j$ with $\Var(\ba r_j) = \frac{\Var(\mu_1) + \Var(\sigma_1)q_j^2 + \Var(\sigma_1 \varepsilon_{1,j})}{n}$. On the other hand, the primitive rating estimator adjusts the scale per user, reporting $q_j + \frac{1}{n} \sum_{i=1}^n \varepsilon_{i,j}$, up to error in estimating the empirical CDFs $F_i$. So the primitive rating estimator has bias on the order of $1/\sqrt M$ and variance on the order of $\frac{\Var(\varepsilon_{1,j}) + 1/M}{n}$ per item. In this model, the primitive rating estimator will outperform the average unless the number of raters $n$ is significantly greater than the number of items each person rates and there is not much heterogeneity in user rating scales.

Figure~\ref{simulation-plots}(b) plots the $L^2$ loss in this model. We compare only the average and the primitive rating estimator, since the full rating estimator is designed to output scores on a different scale. We simulated 500 quality scores $q_j \sim N(0.5,0.15^2)$ with noise $\varepsilon_{i,j} \sim N(0,0.05^2)$, shifts $\mu_i \sim N(0,0.2^2)$, scaling $\sigma_i \sim \text{Unif}(0.25,1.75$) and removed ratings iid with 95\% probability. When users rate 25 items on average, the number of raters needs to be about 3500 (140 times that value) for the average to break even.

\paragraph*{Acknowledgements.} The author gratefully acknowledges Elsa Cazelles, Sinho Chewi, Persi Diaconis, Peng Ding, Steve Evans, Shirshendu Ganguly, Ryan Giordano, Aditya Guntuboyina, Frederick Law, Vilas Winstein, and two anonymous reviewers for many helpful comments. While writing this paper, the author was supported by a Two Sigma PhD Fellowship and a research contract with Sandia National Laboratories, a U.S.\ Department of Energy multimission laboratory. 

\begin{appendix}
\section{Proof of Lemma 3.3 and Proposition 3.5} \label{wasserstein_background_appendix}

\begin{customlem}{3.3} \label{increasing-transport-formula}
Let $\mu,\nu$ be measures on $\R$, and let $f : \R \to \R$ be a function such that $\nu = f_*\mu$. If $f$ is strictly increasing $\mu$-a.e., then $f = F_\nu^{-1} \circ F_\mu$ $\mu$-a.e., where $F_\mu,F_\nu$ are the respective cumulative distribution functions of $\mu,\nu$.
\end{customlem}

\begin{proof}
Let $Y \sim \nu$, so that $f^{-1}(Y) \sim \mu$, where $f^{-1}$ is well-defined $\nu$-a.e. Then
\begin{align*}
F_\mu(x) &= \P( f^{-1}(Y) \leq x) \\
&= \P(Y \leq f(x)) \\
&= F_\nu \circ f(x).
\end{align*}
Therefore, $F_\mu = F_\nu \circ f$ for $\mu$-a.e.\ $x \in [0,1]$. Applying $F_{\nu}^{-1}$ to both sides gives $f = F_\nu^{-1} \circ F_\mu$.
\end{proof}

\begin{customprop}{3.5}[$W_2$-Fr\'echet mean formula] \label{w2-frechet-mean-formula}
Let $\mu_1,\dots,\mu_n$ be measures on $[0,1]$. The empirical Fr\'echet mean of $\mu_1,\dots,\mu_k$ with respect to the Wasserstein 2-distance is the measure $\wh \mu$ with inverse CDF
$$F_{\wh \mu}^{-1} = \frac{1}{n} \sum_{k=1}^n F_k^{-1}.$$
\end{customprop}

\begin{proof}
First, we embed the Wasserstein space $\mc P_2([0,1])$ into $L^2([0,1])$ as follows. Consider the map $Q : \mc P_2([0,1]) \to L^2([0,1],\text{Leb})$ sending $\mu \mapsto F_\mu^{-1}$. Lemma 3.2 tells us that $Q$ is an isometry, so it suffices to find the Fr\'echet mean of $F_1^{-1},\dots, F_n^{-1}$ in $L^2$. So we want to find the $L^2$ function $g$ minimizing $\frac{1}{n} \sum_{k=1}^n \| F_k ^{-1}- g \|_{L^2}^2$. 

The standard bias-variance decomposition calculation tells us that
$$\frac{1}{n} \sum_{k=1}^n \| F_k ^{-1}- g \|_{L^2}^2 = \frac{1}{n} \sum_{k=1}^n \norm{F_k^{-1} - \frac{1}{n} \sum_{j=1}^n F_j^{-1}}_{L^2}^2 + \norm{g - \frac{1}{n} \sum_{j=1}^n F_j^{-1} }_{L^2}^2.$$
The right hand side is uniquely minimized when $g = \frac{1}{n} \sum_{j=1}^n F_j^{-1}$, so this must be the Fr\'echet mean.
\end{proof}

\section{Proof of Proposition 4.2} \label{rating_estimation_with_universal_preferences_appendix}

\begin{customprop}{4.2}[Consistency of the average and primitive rating estimators with $L^2$-loss and universal preferences] \label{R0_equals_A_and_is_consistent}
Let $\mu$ be a measure on $[0,1]$, and let $\phi_1,\phi_2,\ldots : [0,1] \to [0,1]$ be iid random functions. Let $\mu_k := (\phi_k)_*\mu$ for each $1\leq k \leq n$, and denote the cumulative distribution functions of $\mu_k,\mu$ as $F_k,F_\mu$, respectively. Define the assumptions
\begin{itemize}

\item[] (A1). (Consensus preference order): $\mu$ is the population $W_2$-barycenter of the law of $\mu_1$.

\item[] (A2''). (Universal preferences): $\P$-a.s., $\phi_k$ is increasing $\mu$-a.e.

\end{itemize}
Define the rating distribution estimator $\wh \mu$ as the Fr\'echet mean of the $\mu_k$ with respect to the Wasserstein 2-distance; i.e.\ $F_{\wh \mu}^{-1} = \frac{1}{n} \sum_{k=1}^n F_k^{-1}$, where $F_{\wh \mu}$ and $F_{k}$ are the cumulative distribution functions of $\wh \mu, \mu_k$, respectively. Define the average estimator
$$A := \frac{1}{n} \sum_{k=1}^n \phi_k$$
and the primitive rating estimator
$$R_0 := \frac{1}{n} \sum_{k=1}^n F_{\wh \mu}^{-1} \circ F_k \circ \phi_k.$$
Then
\begin{enumerate}[label = (\alph*)]

\item Under assumption (A2''), $A = R_0 = F_{\wh \mu}^{-1} \circ F_\mu$ $\mu$-a.e.,

\item $\| F_{\wh \mu}^{-1} \circ F_\mu - \id \|_{L^2(\mu)} = W_2(\wh \mu,\mu)$.
\end{enumerate}
Consequently, under assumptions (A1), (A2''),
$$\| A - \id \|_{L^2(\mu(m))} = \| R_0 - \id \|_{L^2(\mu(m))} \xrightarrow{\P\text{-a.s.},L^p(\P)} 0 \qquad \text{for all} \quad 1 \leq p < \infty$$
as $n \to \infty$.
\end{customprop}

\begin{proof}
The key to both of these properties is that by Lemma~\ref{increasing-transport-formula}, $\phi_k = F_k^{-1} \circ F_\mu$ $\mu$-a.e.
\begin{enumerate}[label = (\alph*)]

\item Write
\begin{align*}
R_0 &= \frac{1}{n} \sum_{k=1}^n F_{\wh \mu}^{-1} \circ F_k \circ (F_k^{-1} \circ F_\mu) \\
&= \frac{1}{n} \sum_{k=1}^n F_{\wh \mu}^{-1}  \circ F_\mu \\
&= F_{\wh \mu}^{-1} \circ F_\mu,
\end{align*}
\begin{align*}
A &= \frac{1}{n} \sum_{k=1}^n (F_k^{-1} \circ F_\mu) \\
&= \paren{\frac{1}{n} \sum_{k=1}^n F_k^{-1}} \circ F_\mu \\
&= F_{\wh \mu}^{-1} \circ F_\mu.
\end{align*}

\item $F_{\wh \mu}^{-1} \circ F_\mu$ is the optimal transport map from $\mu$ to $\wh \mu$, so
\begin{equation*}
\| F_{\wh \mu}^{-1} \circ F_\mu - \id \|_{L^2(\mu)} = \sqrt{\int_0^1 (F_{\wh \mu}^{-1} \circ F_\mu(x) - x)^2 \, d\mu(x)} = W_2(\wh \mu,\mu).
\end{equation*}

\end{enumerate}
Convergence follows from Lemma 5.
\end{proof}

Statement (b) can alternatively be proven by making a change of variables and applying the CDF formula for Wasserstein distance, Lemma 1.

\section{Calculations for Example 4}

Moving on to the general case where users can disagree on matters of both preference and scale, the following simple example shows us that $A,R_0$ can differ and that neither of them is a consistent estimator of $\id$ in general.

\begin{customex}{4} \label{scaling_and_reversing}
The key feature of this example is that the personal scale functions $\phi_k$ will only be comprised of rescaling and reversing preferences. Let $\mu$ be the uniform distribution on $[1/4,3/4]$, and let $\phi_k(x) = \alpha_k(x-1/2) + 1/2$, where the $\alpha_k$ are picked iid from some distribution with $|\alpha_k| \leq 2$ almost surely.

It follows that the optimal transport maps pushing forward $\mu$ to $\mu_k$ are $\phi_k^*(x) = |\alpha_k|(x-1/2) + 1/2$, and therefore the optimal transport map pushing forward $\mu$ to $\wh \mu$ is $F_{\wh \mu}^{-1} \circ F_\mu(x) = (\frac{1}{n} \sum_{k=1}^n |\alpha_k|)(x-1/2) + 1/2$. We also know
$$F_{\wh \mu}^{-1} \circ F_k(x) = F_{\wh \mu}^{-1} \circ F_\mu \circ (F_k^{-1} \circ F_\mu)^{-1}(x) = \Bigg(\frac{1}{n} \sum_{j=1}^n |\alpha_j| \Bigg) \frac{1}{|\alpha_k|}\paren{x-\frac{1}{2}} + \frac{1}{2}.$$
We can directly calculate the relevant quantities of interest:
\begin{align*}
W_2(\wh \mu,\mu) &= \|F_{\wh \mu}^{-1} \circ F_\mu - \id \|_{L^2(\mu)} \\
&= \sqrt{\int_{1/4}^{3/4} \paren{\frac{1}{n} \sum_{k=1}^n |\alpha_k|(x-1/2) + 1/2 - x}^2 \, dx} \\
&= \sqrt{\int_{-1/4}^{1/4} \paren{\frac{1}{n} \sum_{k=1}^n |\alpha_k|y -y}^2 \, dy} \\
&= \abs{\frac{1}{n} \sum_{k=1}^n |\alpha_k| - 1} \sqrt{\int_{-1/4}^{1/4} y^2 \, dy} \\
&= \frac{1}{4 \sqrt 6} \abs{\frac{1}{n} \sum_{k=1}^n |\alpha_k| - 1},
\end{align*}
We know by assumption (A1) and Lemma 4.1 that $W_2(\wh \mu, \mu) \to 0$, so the strong law of large numbers tells us that $\E[|\alpha_1|] = 1$. Now compare that quantity to the loss for the estimators $A$ and $R_0$:

\begin{align*}
\| A - \id \|_{L^2(\mu)} &= \sqrt{\int_{1/4}^{3/4} \paren{\frac{1}{n} \sum_{k=1}^n \alpha_k(x-1/2) + 1/2 - x}^2 \, dx} \\
&= \sqrt{\int_{-1/4}^{1/4} \paren{\frac{1}{n} \sum_{k=1}^n \alpha_ky - y}^2 \, dy} \\
&= \abs{\frac{1}{n} \sum_{k=1}^n \alpha_k - 1} \sqrt{\int_{-1/4}^{1/4} y^2 \, dy} \\
&= \frac{1}{4 \sqrt 6} \abs{\frac{1}{n} \sum_{k=1}^n \alpha_k - 1},
\end{align*}

\begin{align*}
\| R_0 - \id \|_{L^2(\mu)} &= \sqrt{\int_{1/4}^{3/4} \paren{\paren{\frac{1}{n} \sum_{j=1}^n |\alpha_j|}\frac{1}{n} \sum_{k=1}^n \sgn(\alpha_k)(x-1/2) + 1/2 - x}^2 \, dx} \\
&= \sqrt{\int_{-1/4}^{1/4} \paren{\paren{\frac{1}{n} \sum_{j=1}^n |\alpha_j|}\frac{1}{n} \sum_{k=1}^n \sgn(\alpha_k)y- y}^2 \, dy} \\
&= \abs{\paren{\frac{1}{n} \sum_{j=1}^n |\alpha_j|} \paren{\frac{1}{n} \sum_{k=1}^n \sgn(\alpha_k)} -1} \sqrt{\int_{-1/4}^{1/4} y^2 \, dy} \\
&= \frac{1}{4 \sqrt 6} \abs{\paren{\frac{1}{n} \sum_{j=1}^n |\alpha_j|} \paren{\frac{1}{n} \sum_{k=1}^n \sgn(\alpha_k)} -1}.
\end{align*}
If $\P(\alpha_1 < 0) > 0$, then neither $\| A - \id \|_{L^2(\mu)}$ nor $\| R_0 - \id \|_{L^2(\mu)}$ converges to 0 as $n \to \infty$. In other words, if users are allowed to disagree at all in their preferences, then neither $A$ nor $R_0$ will accurately recover the numerical quality scores.
\end{customex}

\section{Derivation of the rating estimator}

We include here a more thorough derivation of the rating estimator, using Brenier's polar factorization theorem. The analysis shows that given desired properties for our estimator, the definition of $R$ is essentially determined.

\begin{remark} \label{derivation_of_rating_estimator}
An estimator $E$ (considered as a function of the consensus numerical quality $m$) that recovers both the scale and the order of the consensus numerical quality scores should ideally have the following two properties:

\begin{enumerate}[label = (\alph*)]

\item (Scale) $E_*\mu = \wh \mu$: An estimator which accurately recovers individual numerical quality scores should also recover the distribution $\mu$ of scores, and assumption (A1) tells us that $\wh \mu$ is the natural estimator to recover $\mu$.

\item (Order) $E$ is increasing: If $E$ does not (at least approximately) preserve the order of the consensus scores, then it is not reliable for comparing the quality of items.

\end{enumerate}
By Lemma~\ref{increasing-transport-formula}, these two properties imply that $E = F_{\wh \mu}^{-1} \circ F_\mu$.

Now suppose that we want a rating estimator of the form $E = \frac{1}{n} \sum_{k=1}^n f_k \circ \phi_k$ for some choice of functions $f_k$, i.e.\ an estimator where we average the ratings given to us by the users, perhaps post-processing the ratings before averaging them. There are two general strategies for making the sum $E$ equal the sum $F_{\wh \mu}^{-1} \circ F_\mu = \frac{1}{n} \sum_{k=1}^n F_k^{-1} \circ F_\mu$:
\begin{enumerate}[label = (\alph*)]

\item Make the sums equal term by term: $f_k \circ \phi_k = F_k^{-1} \circ F_\mu$.

\item Make each term $f_k \circ \phi_k$ equal to $F_{\wh \mu}^{-1} \circ F_\mu$.

\end{enumerate}
To see what $f_k$ needs to be in each of these cases, apply Brenier's polar factorization theorem, Lemma 3.6, to $\phi_k^{-1}$. This tells us that $\phi_k^{-1} = F_\mu^{-1} \circ F_k \circ \tau_k$, where $\tau_k$ preserves $\mu_k$. Inverting both sides tells us that $\phi_k = \tau_k^{-1} \circ F_k^{-1} \circ F_\mu$. Substituting this into our definition of $E$ gives
$$E = \frac{1}{n} \sum_{k=1}^n f_k \circ \tau_k^{-1} \circ F_k^{-1} \circ F_\mu.$$
So using strategy (a), we need $f_k = \tau_k$, and using strategy (b), $f_k = F_{\wh \mu}^{-1} \circ F_k \circ \tau_k$. When $\phi_k$ is increasing (so $\tau_k = \id$), (a) recovers the average $A$, and (b) recovers the primitive rating estimator $R_0$.

To estimate $\tau_k$, observe that $\tau_k = F_k^{-1} \circ F_\mu \circ \phi_k^{-1}$, which we may write as
$$\tau_k = F_k^{-1} \circ F_\mu \circ {\id} \circ \phi_k^{-1}.$$
The statistician does not have access to $F_\mu$, so we could try using $F_{\wh \mu}$. However, we want our estimate $\wh \tau_k$ to preserve $\mu_k$, so we would need to estimate $\id$ by a map which pushes forward $\mu$ to $\wh \mu$. We don't have such a map (and in fact that is the problem that $E$ is supposed to solve!); the most we have is the primitive rating estimator $R_0$, which asymptotically recovers the ordering of the ratings, and $F_{\nu}$, which we can obtain by calculating the primitive rating of all items. Therefore, we are essentially forced into using the estimator $\wh \tau_k := F_k^{-1} \circ F_{\nu} \circ R_0 \circ \phi_k^{-1}$. Plugging in this estimate of $\tau_k$ into the equation $E = \frac{1}{n} \sum_{k=1}^n f_k \circ \phi_k$ (using either strategy (a) or (b)) gives the estimator
$$E = F_{\wh \mu}^{-1} \circ F_{\nu} \circ R_0,$$
which we recognize as our rating estimator $R$.
\end{remark}

\section{Proof of Theorem 5.2 and Theorem 5.7} \label{asymptotic_consistency_subsection}

The goal of this section is to prove the following asymptotic consistency result for the rating estimator.

\begin{customthm}{5.2}[Consistency of the rating estimator with $L^2$ loss] \label{ratings_are_consistent}
Let $\mu$ be a measure on $[0,1]$, and let $\phi_1,\phi_2,\ldots : [0,1] \to [0,1]$ be iid random functions. Let $\mu_k := (\phi_k)_*\mu$ for each $1\leq k \leq n$, and denote the cumulative distribution functions of $\mu_k,\mu$ as $F_k,F_\mu$, respectively. Let  $g(x) := \E[F_\mu^{-1} \circ F_k \circ \phi_k(x)]$. Suppose that
\begin{itemize}

\item[] (A1). (Consensus scale): $\mu$ is the population $W_2$-Fr\'echet mean of the law of $\mu_1$.

\item[] (A2'). (Preference order regularity): The functions $g,F_\mu^{-1} \circ F_k \circ \phi_k$ are differentiable with $g' \geq \alpha > 0$ and $(F_\mu^{-1} \circ F_k \circ \phi_k)' \geq -\beta$.

\end{itemize}
Define the rating distribution estimator $\wh \mu$ as the Fr\'echet mean of the $\mu_k$ with respect to the Wasserstein 2-distance; i.e.\ $F_{\wh \mu}^{-1} = \frac{1}{n} \sum_{k=1}^n F_k^{-1}$, where $F_{\wh \mu}$ and $F_{k}$ are the cumulative distribution functions of $\wh \mu, \mu_k$, respectively. Define the primitive rating estimator
$$R_0 := \frac{1}{n} \sum_{k=1}^n  F_{\wh \mu}^{-1} \circ F_k \circ \phi_k$$
and the rating estimator
$$R := F_{\wh \mu}^{-1} \circ F_{\nu} \circ R_0,$$
where $\nu := (R_0)_*\mu$. Then
$$\| R - \id \|_{L^2(\mu(m))} \xrightarrow{\P\text{-a.s.},L^p(\P)} 0 \qquad \text{for all} \quad 1 \leq p < \infty$$
as $n \to \infty$.

If $\mu$ is purely atomic, then assumption (A2') can be weakened to
\begin{itemize}

\item[] (A2). (Consensus preference order): The function $g(x) := \E[F_{\mu}^{-1} \circ F_1 \circ \phi_1(x)]$ is strictly increasing $\mu$-a.e.

\end{itemize}

\end{customthm}

Before we can prove this theorem, we first need to take care of a loose end. Looking at the function $g$ from Assumption (A2), the strong law of large numbers tells us that $\frac{1}{n} \sum_{k=1}^n F_\mu^{-1} \circ F_k \circ \phi_k \to g$ pointwise, recovering the true ordering of scores as $n \to \infty$, but we don't actually have access to $F_\mu^{-1}$; instead, we use the bootstrap inverse CDF $F_{\wh \mu}^{-1}$. We must take care to ensure that this bootstrap estimate does not incur much error, even when feeding the bootstrap inverse CDF inputs from random functions.

\begin{lemma} \label{bootstrap-inverse-CDF-is-ok}
For each $x \in [0,1]$, the following holds $\P$-a.e.:
$$R_0(x) := \frac{1}{n} \sum_{k=1}^n F_{\wh \mu}^{-1} \circ F_k \circ \phi_k(x) \to g(x)$$
as $n \to \infty$.
\end{lemma}

This result holds $\P$-almost surely for each point, but it does not necessarily hold $\P$-a.s. for all $x \in [0,1]$ simultaneously. However, using the regularity afforded to us by assumption (A2'), we will be able to upgrade this to a more uniform understanding of the accuracy of the ordering $R_0$ recovers.

We will keep the error from being amplified by function composition by establishing \emph{uniform} convergence of $F_{\wh \mu}^{-1}$ to $F_\mu^{-1}$. The key is the following general result about Wasserstein barycenters on $[0,1]$:

\begin{customlem}{5.4}[Glivenko--Cantelli-type theorem for inverse CDFs of Wasserstein barycenters on ${[0,1]}$]
\label{inverse_CDF_GC}
Let $\Lambda$ be a random measure on $[0,1]$, and let $\lambda$ be the population $W_2$-barycenter of the law of $\Lambda$. Draw iid samples $\Lambda_1,\Lambda_2,\dots$ according to the law of $\Lambda$, and let $ \lambda_n$ be the empirical $W_2$-barycenter of $\Lambda_1,\dots,\Lambda_n$. Denote the cumulative distribution functions of $ \lambda_n$ and $\lambda$ as $F_{ \lambda_n}, F_\lambda$, respectively, and define the generalized inverse cumulative distribution functions as
$$F_{ \lambda_n}^{-1}(x) := \inf \{ t : F_{ \lambda_n}(t) \geq x\}, \qquad F_{\lambda}^{-1}(x) := \inf \{ t : F_{\lambda}(t) \geq x\}.$$
Then
$$\sup_{x \in [0,1]} | F_{ \lambda_n}^{-1}(x) - F_{\lambda}^{-1}(x)| \to 0$$
almost surely as $n \to \infty$. Moreover, if $\supp \lambda = [0,1]$, then
$$\E \squa{\sup_{x \in [0,1]} | F_{ \lambda_n}^{-1}(x) - F_{\lambda}^{-1}(x)|} \leq 4 \sqrt{\frac{\log n}{n}}.$$
for all $n \geq 3$.
\end{customlem}

We postpone the proof of Lemma~\ref{inverse_CDF_GC} to Section~\ref{GC-appendix} of this appendix. Assuming Lemma~\ref{inverse_CDF_GC}, the proof of Lemma~\ref{bootstrap-inverse-CDF-is-ok} is quick.

\begin{proof}[Proof of Lemma~\ref{bootstrap-inverse-CDF-is-ok}]
By the strong law of large numbers,
$$\frac{1}{n} \sum_{k=1}^n F_{ \mu}^{-1} \circ F_k \circ \phi_k(x) \to g(x)$$
$\P$-a.s.\ as $n \to \infty$. To show that $R_0(x)$ has the same limit, Lemma~\ref{inverse_CDF_GC} tells us that $F_{\wh \mu}^{-1} \to F_\mu^{-1}$ uniformly as $n \to \infty$. So given $\varepsilon > 0$, let $N$ be large enough such that for $n \geq N$, $\sup_{y \in [0,1]} |F_{\wh \mu}^{-1}(y) - F_\mu^{-1}(y)| \leq \varepsilon$. Then for all such $n$,
\begin{align*}
&\abs{R_0(x) - \frac{1}{n} \sum_{k=1}^n F_{ \mu}^{-1} \circ F_k \circ \phi_k(x)} \\
&\qquad \qquad \leq \frac{1}{n} \sum_{k=1}^n | F_{\wh \mu}^{-1} \circ F_k \circ \phi_k(x) - F_{ \mu}^{-1} \circ F_k \circ \phi_k(x)| \\
&\qquad \qquad \leq \varepsilon.\qedhere
\end{align*}
\end{proof}

Now we are prepared to prove the consistency of the rating estimator $R$:

\begin{proof}[Proof of Theorem~\ref{ratings_are_consistent}]
We need only prove $\P$-a.s.\ convergence because this implies convergence in probability. Then, for any $1 \leq p < \infty$ and for any $\varepsilon > 0$,
\begin{align*}
\E[\| R - \id \|_{L^2(\mu)}^p] \leq \varepsilon^p \cdot \P(\| R - \id \|_{L^2(\mu)}^p \leq \varepsilon) + 1 \cdot \P(\| R - \id \|_{L^2(\mu)}^p > \varepsilon),
\end{align*}
where the right hand side tends to $\varepsilon^p$ as $n \to \infty$. Here, we have used the fact that $\| R - \id \|_{L^2(\mu)}^p \leq 1$. This gives convergence in $L^p(\P)$.

To decouple the roles of scale and order in the estimation, Brenier's polar factorization theorem tells us that there exists a unique $\mu$-preserving map $\sigma : [0,1] \to [0,1]$ such that $R_0 = F_{\nu}^{-1} \circ F_\mu \circ \sigma$. Substituting this decomposition into the definition of $R$ gives us the simplified theoretical representation
$$R = F_{\wh \mu}^{-1} \circ F_\mu \circ \sigma.$$
Now we can analyze the rating estimation error in two terms which represent the scaling error and ordering error, respectively.
\begin{align*}
\| R - \id \|_{L^2(\mu)} &\leq \| R - \sigma \|_{L^2(\mu)} + \| \sigma - \id \|_{L^2(\mu)} \\
&= \| F_{\wh \mu}^{-1} \circ F_\mu \circ \sigma - \sigma \|_{L^2(\mu)} + \| \sigma - \id \|_{L^2(\mu)} \\
\shortintertext{Using the fact that $\sigma$ is $\mu$-preserving,}
&= \| F_{\wh \mu}^{-1} \circ F_\mu - \id \|_{L^2(\mu)} + \| \sigma - \id \|_{L^2(\mu)}
\shortintertext{By Proposition~\ref{R0_equals_A_and_is_consistent}, the former term is just $W_2(\wh \mu,\mu)$.}
&= W_2(\wh \mu,\mu) + \| \sigma - \id \|_{L^2(\mu)}.
\end{align*}
The first term goes to 0 $\P$-a.s.\ by assumption (A1) and Lemma 4.1, so we now need only concern ourselves with the ordering error term.

To deal with $\| \sigma - \id \|_{L^2(\mu)}$, we first address the case of a purely atomic $\mu$ under assumption (A2). The key point is that the rearrangement $\sigma$ must map atoms of $\mu$ to atoms of $\mu$, so if we preserve the ordering of the atoms, they must be fixed in place. The measure $\mu$ can have at most countably infinitely many atoms, so given $\varepsilon > 0$, let $A_\varepsilon$ be a \emph{finite} set of atoms of $\mu$ such that the total mass of atoms not in $A_\varepsilon$ is $\leq \varepsilon$; choosing $A_\varepsilon$ so that $\mu(\{x \}) > \mu(\{ y \})$ for all $x \in A_\varepsilon$ and $y \notin A_\varepsilon$, we must have $\sigma(A_\varepsilon) = A_\varepsilon$. Then, letting $\delta = \min \{ |g(x)-g(y)| : x,y \in A_\varepsilon \}$, we look to Lemma~\ref{bootstrap-inverse-CDF-is-ok} to guarantee that $\P$-a.e., for all sufficiently large $n$, $|R_0(x) - g(x)|  < \delta/2$ for all $x \in A_\varepsilon$. Since $g$ is increasing, the ordering of the $x \in A_\varepsilon$ is the same as the ordering of the $g(x)$ in $g(A_\varepsilon)$. Moreover, the choice of $\delta$ guarantees that the $R_0(x)$ preserve the same ordering because if $x < y$, then
$$R_0(y) - R_0(x) > g(y) - g(x) - 2\delta \geq 0.$$
Finally, recalling that $R_0 = F_{\wt \mu}^{-1} \circ F_\mu \circ \sigma$, where $F_{\wt \mu}^{-1} \circ F_\mu$ is nondecreasing, we see that $\sigma$ preserves the ordering of atoms in $A_\varepsilon$ for all sufficiently large $n$, as well. The conclusion is that since $\sigma$ preserves the ordering of the points in the finite set $A_\varepsilon$ and maps $A_\varepsilon$ to $A_\varepsilon$ bijectively, $\sigma$ must act as the identity on $A_\varepsilon$ for all sufficiently large $n$, $\P$-a.s. This proves the statement in the case where $\mu$ is purely atomic. 

Assuming (A2'), given an arbitrary $\varepsilon > 0$, we will show that with high probability, for all sufficiently large $n$, $\sigma$ moves points at most $\varepsilon$. The strategy is to have $R_0$ preserve the ordering of a fine mesh of points in $[0,1]$ and then use the regularity assumption (A2) to extend this to preserving order between most pairs of points; this preservation of order will be captured by $\sigma$, which will then be close to the identity map.  Denote $R_1 := \frac{1}{n} \sum_{k=1}^n F_\mu^{-1} \circ F_k \circ \phi_k$, and let $m$ be large enough such that $3\beta/(\alpha m) + 2/m \leq \varepsilon$. The strong law of large numbers tells us that for all sufficiently large $n$, $|R_1(j/m) - g(j/m)| < \delta$ for all $0 \leq j \leq m$, where we choose $\delta = \beta/(8m)$. We may also let $n$ be large enough such that $\sup_{x \in [0,1]} |R_0(x) - R_1(x)| < \delta$, by Lemma~\ref{inverse_CDF_GC}.

Now let $x,y \in [0,1)$ with $x < y$, and denote $k, \ell$ as the indices such that $(k-1)/M \leq x < k/M$ and $(\ell-1)/M \leq y < \ell/M$. Then if $|x-y| > \varepsilon$ (so that $(\ell-1 -k)/M \geq 3\beta/(\alpha M)$),
\begin{align*}
&R_0(y) - R_0(x) \\
& \qquad \geq R_1(y) - R_1(x) - 2\delta \\
& \qquad = \squa{R_1(\tfrac{\ell-1}{m}) + (R_1(y) - R_1(\tfrac{\ell-1}{m}))} - \squa{R_1(\tfrac{k}{m}) - (R_1(\tfrac{k}{m}) - R_1(x))} - 2\delta \\
& \qquad = \squa{R_1(\tfrac{\ell-1}{m}) + \int_{\frac{\ell-1}{m}}^y R_1'(s) \, ds} - \squa{R_1(\tfrac{k}{m}) - \int_x^{\frac{k}{m}} R_1'(s) \, ds} - 2\delta 
\intertext{Since each $(F_\mu^{-1} \circ F_k \circ \phi_k)' \geq - \beta$, $R_1' \geq -\beta$, as well.}
& \qquad \geq R_1(\tfrac{\ell-1}{m}) - R_1(\tfrac{k}{m}) - \beta(y - \tfrac{\ell-1}{m}) - \beta(\tfrac{k}{m} - x) - 2\delta \\
& \qquad \geq g(\tfrac{\ell-1}{m}) - g(\tfrac{k}{m}) - \beta(y - \tfrac{\ell-1}{m}) - \beta(\tfrac{k}{m} - x) - 4\delta \\
& \qquad \geq \alpha(\tfrac{\ell-1-k}{m}) - \beta(y - \tfrac{\ell-1}{m}) - \beta(\tfrac{k}{m} - x) - 4\delta \\
& \qquad \geq \frac{3\beta}{m} - \frac{2\beta}{m} - 4\delta \\
& \qquad = \frac{\beta}{2m} \\
& \qquad > 0.
\end{align*}
Keep in mind that this conclusion holds $\P$-a.s.\ for all such pairs $x,y$ simultaneously.

We have proven that for any $\varepsilon$, there exists a sample size $N$ such that for all $n \geq N$,
$$|x-y| > \varepsilon \text{ and } x < y \implies R_0(y) > R_0(x) \implies \sigma(y) > \sigma(x).$$
We then have that for $x \in (\varepsilon, 1-\varepsilon)$, $\sigma(x) > \sigma(y)$ for all $y \in [0,x-\varepsilon]$ and $\sigma(x) < \sigma(y)$ for all $y \in [x+\varepsilon,1]$. So, since $\sigma$ preserves $\mu$, we get
$$F_\mu^{-1} \circ F_\mu(x - \varepsilon) \leq \sigma(x) \leq F_\mu^{-1} \circ F_\mu(x + \varepsilon),$$
Since $\sigma$ maps $\supp \mu$ to itself, we may conclude that that for $\mu$-a.e.\ $x \in [0,1]$,
$$x - \varepsilon \leq \sigma(x) \leq x + \varepsilon.$$
So for all sufficiently large $n$,
$$\| \sigma - \id\|_{L^2(\mu)} \leq\varepsilon$$
$\P$-a.s. Since $\varepsilon >0$ was arbitrary, we get that $\| \sigma - \id\|_{L^2(\mu)} \to 0$ $\P$-a.s.\ as $n \to \infty$, as desired.
\end{proof}

\begin{customthm}{5.7}[Rates of convergence for the rating estimator] \label{rating_rates_of_convergence}
With the same notation as Theorem~\ref{ratings_are_consistent},
\begin{enumerate}[label = (\alph*)]

\item If $\mu$ is purely atomic with $M$ atoms, assuming (A1) and (A2),
\begin{align*}
\E[\| R - \id \|_{L^2(\mu)}] &\leq \frac{1}{2\sqrt n} + 4M e^{-n\delta^2/8} \qquad \forall n \geq 1 \\
& = O \paren{\frac{1}{\sqrt n}},
\end{align*}
where $\delta := \min \{ g(y) - g(x) : x,y \text{ are atoms of $\mu$}, x < y  \}$.

\item If $\supp \mu = [0,1]$, assuming (A1) and (A2'),
\begin{align*}
\E[\| R - \id \|_{L^2(\mu)}] &\leq \paren{\frac{3C}{\alpha} + \frac{11}{2}} \sqrt{\frac{\log n}{n}} \qquad \forall n \geq 3 \\
&= O \paren{\sqrt{\frac{\log n}{n}}},
\end{align*}
where $C$ is a constant depending only on $\beta$.
\end{enumerate}
\end{customthm}

\begin{proof}
Writing $\E[\| R - \id \|_{L^2(\mu)}] \leq \E[W_2(\wh \mu,\mu)] + \E[\| \sigma - \id \|_{L^2(\mu)}]$, we bound each term. The first term can be bounded by
\begin{align*}
\E[W_2(\wh \mu,\mu)] &= \E[ \|F_{\wh \mu}^{-1} - F_\mu^{-1}\|_{L^2(dx)}] \\
&= \E \squa{\norm{\frac{1}{n} \sum_{k=1}^n F_k^{-1} - F_\mu^{-1}}_{L^2(dx)}} \\
&\leq \sqrt{\E \squa{\norm{\frac{1}{n} \sum_{k=1}^n F_k^{-1} - F_\mu^{-1}}_{L^2(dx)}^2}} \\
&= \sqrt{\E \squa{\frac{1}{n^2} \sum_{j,k=1}^n \int_0^1 (F_k^{-1}(x) - F_\mu^{-1}(x))(F_j^{-1}(x) - F_\mu^{-1}(x)) \, dx}} \\
&= \sqrt{\frac{1}{n^2} \sum_{j,k=1}^n \int_0^1 \E[ (F_k^{-1}(x) - F_\mu^{-1}(x))(F_j^{-1}(x) - F_\mu^{-1}(x))] \, dx} \\
\intertext{Assumption (A1) implies that $\E[F_1^{-1}(x)] = F_\mu^{-1}(x)$ for all $x \in [0,1]$ (see Lemma~\ref{pointwise_quantile_convergence}), so noting that the $F_k^{-1}(x) - F_\mu^{-1}(x)$ are centered, iid random variables for each $x$, we see that the off-diagonal terms will all be 0.}
&= \sqrt{\frac{1}{n^2} \sum_{k=1}^n \int_0^1 \Var(F_k^{-1}(x)) \, dx} \\
&\leq \frac{1}{2\sqrt n}.
\end{align*}

For the remaining term, if $\mu$ is purely atomic with $N$ atoms and we are assuming (A2), then we will have $\sigma = \id$ once $|R_0(x) - g(x)| < \delta/2$ for each of the $N$ atoms $x$ of $\mu$, where $\delta$ is a constant depending on $g$. To compare $R_0$ to $g$, it will be convenient to modify $R_0$ by removing the doubled dependence on each $F_k$ in each term. In each term, we replace $F_{\wh \mu}^{-1} := \frac{1}{n} \sum_{j=1}^n F_j^{-1}$ with $F_{\nu_k}^{-1} := \frac{1}{n} \sum_{j \neq k} F_j^{-1} + \frac{1}{n} F_{k'}^{-1}$, where $F_{k'}^{-1}$ is an independent copy of $F_k^{-1}$. This gives $R_2 := \frac{1}{n} \sum_{k=1}^n F_{\nu_k}^{-1} \circ F_k \circ \phi_k$, where $|R_0 - R_2| \leq \frac{1}{n}$. Moreover, $\E[R_2(x)] = g(x)$ because
\begin{align*}
&\E[F_{\nu_k}^{-1} \circ F_k \circ \phi_k(x) \mid \phi_k] \\
&\qquad \qquad = \frac{1}{n} \sum_{j \neq k} \E[F_j^{-1} \circ F_k \circ \phi_k(x) \mid \phi_k] + \frac{1}{n} \E[ F_{k'}^{-1} \circ F_k \circ \phi_k(x) \mid \phi_k] \\
&\qquad \qquad = F_\mu^{-1} \circ F_k \circ \phi_k(x),
\end{align*}
which gives
$$\E[R_2(x)] = \frac{1}{n} \sum_{k=1}^n \E[\E[F_{\nu_k}^{-1} \circ F_k \circ \phi_k(x) \mid \phi_k]] = g(x).$$

Using the bounded differences inequality (see, e.g., Corollary 2.21 of \cite{wainwright2019high}) with bounded difference $\frac{2}{n}$, we have
$$\P(|R_2(x) - g(x)| \geq t) \leq 2 e^{-nt^2/2},$$
so that
\begin{align*}
\P(|R_0(x) - g(x)| \geq \delta/2) &\leq \P(|R_2(x) - g(x)| \geq \delta/2 - 1/n) \\
&\leq 2 e^{-n(\delta/2-1/n)^2/2} \\
&\leq 4 e^{-n\delta^2/8},
\end{align*}
where we have used $e^{\delta^2/2} \leq e^{1/2} \leq 2.$
Applying a union bound over the $M$ atoms, we get
$$\P(|R_0(x) - g(x)| \geq \delta/2 \text{ for any atom $x$}) \leq 4M e^{-n\delta^2/8}.$$
Then
\begin{align*}
\E[\|\sigma - \id \|_{L^2(\mu)}] &\leq 1 \cdot \P(|R_0(x) - g(x)| \geq \delta/2 \text{ for some atom $x$}) \\
&\leq 4M e^{-n\delta^2/8}
\end{align*}
finishing the proof for this case.

For the remaining case, assuming (A2') and $\supp \mu = [0,1]$, we will show that on the event where $\sigma$ approximates $\id$ as in the proof of Theorem~\ref{ratings_are_consistent}, $\| \sigma - \id \|_{L^2(\mu)} \leq (3\beta/\alpha + 2) \sqrt{\frac{\log n}{n}}$. We will also show that the probability of the complement event is appropriately small. Denoting $R_1 := \frac{1}{n} \sum_{k=1}^n F_\mu^{-1} \circ F_k \circ \phi_k$, for the ``good'' event to occur, we need $|R_1(x) - g(x)| \leq \delta :=  \beta/(8m)$ for $m$ values of $x$ and $\sup_{x \in [0,1]} |R_0(x) - R_1(x)| \leq \delta := \beta/(8m)$, where we choose $m = \sqrt{\frac{n}{\log n}}$. Again by Hoeffding's inequality,
$$\P(|R_1(x) - g(x)| \geq \delta) \leq 2 e^{-2n\delta^2},$$
so that
\begin{align*}
&\P(|R_1(x) - g(x)| \geq \delta \text{ for any of the $m$ points $x$}) \\
& \qquad \leq 2 \sqrt{\frac{n}{\log n}} e^{-2n\delta^2} \\
& \qquad = 2 \sqrt{\frac{n}{\log n}} \exp \paren{- \beta^2 \log n/32} \\
& \qquad = \frac{2}{n^{\beta^2/32 - 1/2}\sqrt{\log n}}.
\end{align*}
By Lemma~\ref{inverse_CDF_GC} and the bounded differences inequality,
\begin{align*}
&\P \paren{\sup_{x\in [0,1]} |R_0(x) - R_1(x)| \geq \delta} \\
&\qquad \leq \P \paren{\sup_{x\in [0,1]} |R_0(x) - R_1(x)| - \E \squa{\sup_{x\in [0,1]} |R_0(x) - R_1(x)|}  \geq \delta - 4 \sqrt{\frac{\log n}{n}}} \\
&\qquad \leq \exp \paren{- \frac{n(\delta - 4 \sqrt{\frac{\log n}{n}})^2}{4}} \\
&\qquad = \exp \paren{- \frac{(\beta/8 - 4)^2\log n}{4}} \\
&\qquad = \frac{1}{n^{(\beta/16 - 2)^2}}.
\end{align*}
In total, we get
\begin{align*}
\E[\|\sigma - \id \|_{L^2(\mu)}] &\leq (3\beta/\alpha + 2) \sqrt{\frac{\log n}{n}} \cdot 1 \\
& \qquad \qquad + 1 \cdot \P(|R_1(x) - g(x)| \geq \delta \text{ for any of $M$ points}) \\
& \qquad \qquad + 1 \cdot \P \paren{\sup_{x\in [0,1]} |R_0(x) - R_1(x)| \geq \delta} \\
&= (3\beta/\alpha + 2) \sqrt{\frac{\log n}{n}} + \frac{2}{n^{\beta^2/32 - 1/2}\sqrt{\log n}} + \frac{1}{n^{(\beta/16 - 2)^2}}
\shortintertext{If assumption (A2') holds for $(\alpha,\beta)$, then it holds for $(\alpha,\beta')$ for any $\beta' > \beta$, so if $\beta^2 /32 - 1/2 < 1/2$ or $(\beta / 16 - 2)^2 < 1/2$, we may repeat the argument with the increased value $\beta' = 4 \sqrt 2$ (if $\beta' < 4 \sqrt 2$) or $\beta' = 32 + 8 \sqrt 2$. So we get}
&\leq \paren{\frac{3C}{\alpha} + 5} \sqrt{\frac{\log n}{n}},
\end{align*}
where
$$C = \begin{cases}
4 \sqrt 2 &\text{if } \beta \leq 4 \sqrt 2 \\
\beta &\text{if } 4 \sqrt 2 < \beta \leq 32 - 8 \sqrt 2 \\
32 + 8 \sqrt 2 &\text{if } 32 - 8 \sqrt 2 < \beta \leq 32 + 8 \sqrt 2 \\
\beta &\text{if } \beta > 32 + 8 \sqrt 2.
\end{cases}$$
Adding this to the $\frac{1}{2\sqrt n}$ contribution from $W_2(\wh \mu,\mu)$ completes the proof.
\end{proof}

\section{Proof of Theorem 6.1} \label{incomplete-data-appendix}

\begin{customthm}{6.1} \label{incomplete-consistency-thm}
Let $\mu$ be purely atomic with $M$ atoms, and assume that

\begin{itemize}

\item[] (A1). (Consensus scale): $\mu$ is the population $W_2$-Fr\'echet mean of the law of $\wt \mu_1$.

\item[] ({{A2'}'}'). (Consensus preference order with incomplete data): The function $h(x) := \E[F_\mu^{-1} \circ \wt F_1 \circ \phi_1(x)]$ is increasing.

\item[] (A3). (Balanced choices): $\E[\wt F_1^{-1}(x) \mid \phi_1] = F_1^{-1}(x)$ for all $x \in [0,1]$.

\item[] (A4). (Item choices don't influence profiles): The collection $N := \{ N_x : x \in [0,1]\}$ of ``who rates what'' is independent of $\{ \phi_k : k \geq 1\}$.

\end{itemize}
Then
$$\E[\| \wt R - \id \|_{L^2(\mu)}] \leq \frac{1}{2\sqrt n} + 6.5 \sum_{\text{atoms $x$}} \E[e^{-|N_x|\delta^2/128}],$$
where $\delta := \min \{ h(y) - h(x) : \text{$x,y$ are atoms of $\mu$, $x < y$}\}$. Consequently, if for all atoms $x$ of $\mu$, $|N_x| \geq \frac{64}{\delta^2} \log n$, then
\begin{align*}
\E[\| \wt R - \id \|_{L^2(\mu)}] &\leq \frac{7M}{\sqrt n} \\
&= O\paren{\frac{1}{\sqrt n}}.
\end{align*}
\end{customthm}

\begin{lemma} \label{wasserstein_lln}
$$\E[W_2(\wh \mu,\mu)] \leq \frac{1}{2\sqrt n}.$$
\end{lemma}

We will not prove this lemma, as the proof is a simplified version of the proof of the following lemma, which we will also need.

\begin{lemma} \label{tilde_approx_hat}
Assume that $\E[\wt F_1^{-1}(x) \mid \phi_1] = F_\mu^{-1}(x)$ for all $x \in [0,1]$. Then
$$\E[W_2(\wt \mu,\wh \mu)] \leq \frac{1}{2\sqrt n}.$$
\end{lemma}

\begin{proof}[Proof of Lemma~\ref{tilde_approx_hat}]
\begin{align*}
\E[W_2(\wt \mu,\wh \mu)] &\leq \sqrt{\E[W_2(\wt \mu,\wh \mu)^2]} \\
&= \sqrt{\int_0^1 \E[(F_{\wt \mu}^{-1}(x) - F_{\wh \mu}^{-1}(x))^2] \, dx} \\
&= \sqrt{\frac{1}{n^2} \sum_{j=1}^n \sum_{k=1}^n \int_0^1 \E[(\wt F_j^{-1}(x) - F_{j}^{-1}(x))(\wt F_k^{-1}(x) - F_k^{-1}(x))] \, dx}
\intertext{For $j \neq k$, the two terms in the expectation are independent. Moreover, they are centered, as $\E[\wt F_k^{-1}(x)] = \E[\E[\wt F_k^{-1}(x) \mid \phi_k]] = \E[F_k^{-1}(x)]$. So only the diagonal terms survive.}
&= \sqrt{\frac{1}{n} \int_0^1 \E[(\wt F_1^{-1}(x) - F_1^{-1}(x))^2] \, dx} \\
&= \sqrt{\frac{1}{n} \int_0^1 \E[\E[(\wt F_1^{-1}(x) - F_1^{-1}(x))^2 \mid \phi_1]] \, dx} \\
&= \sqrt{\frac{1}{n} \int_0^1 \E[\Var(\wt F_k^{-1}(x) \mid \phi_1)] \, dx} \\
&\leq \frac{1}{2\sqrt n}.\qedhere
\end{align*}
\end{proof}

\begin{lemma}[Concentration of primitive ratings] \label{concentration_of_primitive_ratings} \label{concentration-incomplete-lemma}
Assume (A2) and (A3). Then for any $x \in [0,1]$ and $0 < t \leq 1$,
$$\P(|\wt R_0(x) - h(x)| \geq t) \leq 6.5 \E[e^{-|N_x| t^2/32}].$$
\end{lemma}

\begin{proof}[Proof of Lemma~\ref{concentration-incomplete-lemma}]
Write
\begin{align*}
|\wt R_0(x) - h(x)| &\leq \Bigg | \underbrace{\underbrace{\frac{1}{|N_x|} \sum_{k \in N_x} F_{\wt \mu}^{-1} \circ \wt F_k \circ \phi_k(x)}_{a} - \underbrace{\frac{1}{n} \sum_{k=1}^n F_{\wt \mu}^{-1} \circ \wt F_k \circ \phi_k(x)}_b}_{I} \Bigg | \\
&\qquad \qquad+ \Bigg | \underbrace{\frac{1}{n} \sum_{k=1}^n F_{\wt \mu}^{-1} \circ \wt F_k \circ \phi_k(x) - h(x)}_{II} \Bigg |.
\end{align*}
To deal with $I$, we will show that each of the two averages is similar to its mean, conditioned on all $N_x$; these conditional expectations are the same, so this will show that the two averages are very similar to each other. Letting $N := \{ N_x : x \in [0,1]\}$, for any $t > 0$,
\begin{align*}
\P(| I | \geq t/2 \mid N) &\leq \P(|a - \E[a \mid N]|  \geq t/4 \mid N) + \P(|b - \E[b \mid N]| \geq t/4 \mid N) \\
&= \E[\P(|a - \E[a \mid N]|  \geq t/4 \mid N, \{ \phi_j : j \notin N_a\}) \mid N] \\
&\qquad + \P(|b - \E[b \mid N]| \geq t/4 \mid N) \\
\shortintertext{Using the bounded differences inequality on each term (with bounded differences $\frac{2}{|N_x|}$ and $\frac{2}{n}$, respectively),}
&\leq \E[2 e^{-|N_x| t^2/32} \mid N] + 2e^{-n t^2/32} \\
&\leq 4 e^{-|N_x|t^2/32}.
\end{align*}
Taking the expectation over $N$, we get
$$\P( | I | \geq t/2) \leq 4 \E[e^{-|N_x| t^2/32}].$$

To deal with $II$, first look at each term $F_{\wt \mu}^{-1} \circ \wt F_k \circ \phi_k(x) = \frac{1}{n} \sum_{\ell=1}^n \wt F_{\ell}^{-1} \circ \wt F_k \circ \phi_k(x)$. We can replace $F_k^{-1}$ by an independent copy $F_{k'}^{-1}$, denoting the resulting term as $F_{\wt \nu_k}^{-1} \circ \wt F_k \circ \phi_k(x)$, at the cost of, at most, $1/n$. That is,
$$|F_{\wt \mu}^{-1} \circ \wt F_k \circ \phi_k(x) - F_{\wt \nu_k}^{-1} \circ \wt F_k \circ \phi_k(x)| \leq \frac{1}{n}.$$
Doing so ensures that $\frac{1}{n} \sum_{k=1}^n F_{\wt \nu_k}^{-1} \circ \wt F_k \circ \phi_k(x) - h(x)$ is centered, as
\begin{align*}
\E[F_{\wt \nu_k}^{-1} \circ \wt F_k \circ \phi_k(x)] &= \E[\E[F_{\wt \nu_k}^{-1} \circ \wt F_k \circ \phi_k(x) \mid \phi_k,\wt F_k]] \\
&= \frac{1}{n} \sum_{\ell \neq k} \E[\E[\wt F_{\ell}^{-1} \circ \wt F_k \circ \phi_k(x) \mid \phi_k,\wt F_k]] \\
&\qquad  + \frac{1}{n} \E[\E[\wt F_{k'}^{-1} \circ \wt F_k \circ \phi_k(x) \mid \phi_k,\wt F_k]] \\
&=\E[ F_\mu^{-1} \circ \wt F_k \circ \phi_k(x)] \\
&= h(x)
\end{align*}
for all $k$. Thus, we may apply the bounded differences inequality (with bounded difference $2/n$), conditioned on $\wt F_{1'}^{-1},\dots,\wt F_{n'}^{-1}$, to get
$$\P\paren{\abs{\frac{1}{n} \sum_{k=1}^n F_{\wt \nu_k}^{-1} \circ \wt F_k \circ \phi_k(x) - h(x)} \geq t/2 \mid N, \wt F_{1'}^{-1},\dots,\wt F_{n'}^{-1}} \leq 2 e^{-nt^2/8}.$$
Note that, as in our previous use of the bounded differences inequality, we are leveraging the fact that $\phi_1, \dots, \phi_n$ are independent given $N$, i.e.\ assumption (A4). Averaging over $\wt F_{1'}^{-1},\dots,\wt F_{n'}^{-1}$ gives
$$\P\paren{\abs{\frac{1}{n} \sum_{k=1}^n F_{\wt \nu_k}^{-1} \circ \wt F_k \circ \phi_k(x) - h(x)} \geq t/2} \leq 2 e^{-nt^2/8},$$
and we replace the $\wt \nu_k$ terms with the original $\wt \mu$ terms as follows:
\begin{align*}
\P(| II| \geq t/2) &= \P\paren{\abs{\frac{1}{n} \sum_{k=1}^n F_{\wt \mu}^{-1} \circ \wt F_k \circ \phi_k(x) - h(x)} \geq t/2} \\
&\leq \P\paren{\abs{\frac{1}{n} \sum_{k=1}^n F_{\wt \mu}^{-1} \circ \wt F_k \circ \phi_k(x) - h(x)} \geq t/2 - 1/n} \\
& \leq 2 e^{-n(t/2-1/n)^2/8}
\end{align*}

Putting the bounds for $I$ and $II$ together, we have
\begin{align*}
\P(| \wt R_0(x) - h(x)| \geq t) &\leq \P(| I | + | II | \geq t) \\
&\leq \P(| I | \geq t/2) + \P(| II | \geq t/2) \\
&\leq 4 \E[e^{-|N_x| t^2/32}] + 2 e^{-n(t/2 -1/n)^2/8} \\
\shortintertext{Using $e^{t/8} \leq e^{1/8} \leq 1.25$,}
&\leq 4 \E[e^{-|N_x| t^2/32}] + 2.5 e^{-nt^2/32} \\
&\leq 6.5 \E[e^{-|N_x| t^2/32}].\qedhere
\end{align*}
\end{proof}

\begin{proof}[Proof of Theorem~\ref{incomplete-consistency-thm}]
Using Brenier's polar factorization theorem, we may write $\wt R_0 := F_{(\wt R_0)_*\mu}^{-1} \circ F_\mu \circ \wt \sigma$, where $\wt \sigma$ preserves $\mu$. Then $\wt R = F_{\wt \mu}^{-1} \circ F_\mu \circ \wt \sigma$, and
\begin{align*}
\| \wt R - \id \|_{L^2(\mu)} &\leq \| F_{\wt \mu}^{-1} \circ F_\mu \circ \wt \sigma - \wt \sigma \|_{L^2(\mu)} + \| \wt \sigma - \id \|_{L^2(\mu)} \\
&= \| F_{\wt \mu}^{-1} \circ F_\mu - \id \|_{L^2(\mu)} + \| \wt \sigma - \id \|_{L^2(\mu)} \\
&\leq \| F_{\wt \mu}^{-1} \circ F_\mu - F_{\wh \mu}^{-1} \circ F_\mu \|_{L^2(\mu)}  + \| F_{\wh \mu}^{-1} \circ F_\mu - \id \|_{L^2(\mu)} + \| \wt \sigma - \id \|_{L^2(\mu)} \\
&= \| F_{\wt \mu}^{-1} - F_{\wh \mu}^{-1} \|_{L^2(\lambda)}  + W_2(\wh \mu,\mu) + \| \wt \sigma - \id \|_{L^2(\mu)} \\
&=W_2(\wt \mu, \wh \mu)  + W_2(\wh \mu,\mu) + \| \wt \sigma - \id \|_{L^2(\mu)}.
\end{align*}
Applying Lemma~\ref{wasserstein_lln} and Lemma~\ref{tilde_approx_hat}, we get
$$\E[\| \wt R - \id \|_{L^2(\mu)}] \leq \frac{1}{\sqrt n} + \E[\| \wt \sigma - \id \|_{L^2(\mu)}].$$

To deal with the remaining term, our strategy is to show that with high probability, $\wt R_0$ and $h$ are about the same. This will tell us that the orderings $\wt R_0$ and $h$ give to the inputs in $[0,1]$ are approximately the same, and, by extension, $\wt \sigma$ approximately preserves the ordering of inputs in $[0,1]$. Since $\wt \sigma$ is $\mu$-preserving, this will force $\wt \sigma \approx \id$ on the support of $\mu$.

Using Lemma~\ref{concentration_of_primitive_ratings} with $t = \delta/2$, we have that for each atom $x$ of $\mu$,
$$\P(| \wt R_0(x) - h(x)| \geq \delta/2) \leq 6.5 \E[e^{-|N_x| \delta^2/128}].$$
A union bound over all $M$ values of $x$ gives
$$\P(| \wt R_0(x) - h(x)| \geq \delta/2 \text{ for any atom $x$}) \leq 6.5 \sum_{\text{atoms } x} \E[e^{-|N_x| \delta^2/128}].$$

Putting all the pieces together gives
\begin{align*}
\E[\| \wt R - \id \|_{L^2(\mu)}] &\leq \frac{1}{2\sqrt n} + \E[\| \wt \sigma - \id \|_{L^2(\mu)}] \\
&\leq \frac{1}{2\sqrt n} + 1 \cdot \P(\wt \sigma \neq \id) \\
&\leq \frac{1}{2\sqrt n} + \P(|\wt R_0(x) - h(x)| \geq \delta/2 \text{ for $\geq 1$ of the $M$ item ratings $x$}) \\
&\leq \frac{1}{2\sqrt n} + 6.5 \sum_{\text{atoms } x} \E[e^{-|N_x| \delta^2/128}].\qedhere
\end{align*}
\end{proof}

\section{Proof of Proposition 7.1 and Proposition 7.2} \label{ranking-section-appendix}

\begin{customprop}{7.1}[Concordance statistic formulas] \label{concordance-formulas}
Given a matrix $(r_{i,j})_{1 \leq i \leq n, 1 \leq j \leq M}$ of rating $i$ of item $j$, let $\ba r_i := \frac{1}{M} \sum_{j=1}^M r_j$ be user $i$'s rating mean, and let $\sigma_i^2 := \frac{1}{M} \sum_{j=1}^M (r_{i,j} - \ba r_i)^2$ be user $i$'s rating variance.
\begin{enumerate}[label = (\alph*)]

\item
\begin{align*}
W_{\on{scale}} &= \frac{\sum_{i=1}^n \sum_{k=1}^n \Cov(X_i,X_k)}{n \sum_{i=1}^n \Var(X_i)} \\
&= \frac{\sum_{i=1}^n \sum_{k=1}^n (\sum_{j=1}^M r_{i,(j)}r_{k,(j)}) - M\ba r_i \ba r_k }{n M \sum_{i=1}^n \sigma_i^2}
\end{align*}
where $X_i \sim \mu_i$ and the $X_i$ are monotone coupled (i.e.\ $X_k = F_k^{-1} \circ F_i(X_i)$ for each $i,k$). Here, $r_{i,(j)}$ means the $j$-th lowest rating from user $i$.

\item
\begin{align*}
W_{\on{ratings}} &= \frac{\sum_{i=1}^n \sum_{k=1}^n \sum_{\ell=1}^n \sum_{p=1}^n \Cov(Y_{i,k},Y_{\ell,p}) }{n^3 \sum_{i=1}^n \Var(X_i)} \\
&= \frac{\sum_{i=1}^n \sum_{k=1}^n \sum_{\ell=1}^n \sum_{p=1}^n (\sum_{j=1}^M r_{k,(\on{ind}(r_{i,j}))} r_{p,(\on{ind}(r_{\ell,j}))}) - M\ba r_i \ba r_k }{n^3 M \sum_{i=1}^n \sigma_i^2},
\end{align*}
where $Y_{i,k} := F_k^{-1} \circ F_i(r_{i,J})$ and $J \sim \on{Uniform}(\{1,\dots,M\})$. Here, $\on{ind}(r_{i,j})$ denotes the index of $r_{i,j}$ among user $i$'s sorted ratings $r_{i,(1)},\dots,r_{i,(M)}$, and $r_{k,(\ind(r_{i,j}))}$ is the $\ind(r_{i,j})$-th smallest rating given by user $k$.

\end{enumerate}
\end{customprop}

\begin{proof}
\begin{enumerate}[label = (\alph*)]

\item[]

\item First, we calculate the expectation of $\wh \mu$:

\begin{align*}
\int x \, d\wh \mu(x) &= \int F_{\wh \mu}^{-1}(y) \, dy \\
&= \frac{1}{n} \sum_{i=1}^n \int F_k^{-1}(y) \, dy \\
&= \frac{1}{n} \sum_{i=1}^n \ba r_i,
\end{align*}
which we denote as $\ba r$. It now remains to calculate $\Var(\wh \mu)$.
\begin{align*}
\int (x - \ba r)^2 \, d\mu(x) &= \int (F_{\wh \mu}^{-1}(y) - \ba r)^2 \, dy \\
&= \int \paren{\frac{1}{n} \sum_{i=1}^n (F_i^{-1}(y) - \ba r_i)}^2 \, dy \\
&= \frac{1}{n^2} \sum_{i=1}^n \sum_{k=1}^n \int (F_i^{-1}(y) - \ba r_i)(F_k^{-1}(y) - \ba r_k) \, dy \\
&= \frac{1}{n^2} \sum_{i=1}^n \sum_{k=1}^n \Cov(X_i,X_k).
\end{align*}

\item Now, we calculate the expectation and variance of $(R_0)_*\mu$:

\begin{align*}
\int R_0(x) \, d\mu(x) &= \frac{1}{n} \sum_{i=1}^n \int F_{\wh \mu}^{-1} \circ F_i \circ \phi_i(x) \, d\mu(x) \\
&= \frac{1}{n} \sum_{i=1}^n \int y \, d\wh \mu(y) \\
&= \ba r.
\end{align*}
\begin{align*}
&\int (R_0(x) - \ba r)^2 \, d\mu(x) \\
&\qquad = \int \paren{\frac{1}{n} \sum_{i=1}^n F_{\wh \mu}^{-1} \circ F_i \circ \phi_i(x) - \ba r }^2 \, d\mu(x) \\
&\qquad = \int \paren{\frac{1}{n^2} \sum_{k=1}^n \sum_{i=1}^n F_k^{-1} \circ F_i \circ \phi_i(x) - \ba r_k }^2 \, d\mu(x) \\
&\qquad = \frac{1}{n^4} \sum_{i=1}^n \sum_{k=1}^n \sum_{\ell=1}^n \sum_{p=1}^n \int (F_k^{-1} \circ F_i \circ \phi_i(x) - \ba r_k)(F_p^{-1} \circ F_\ell \circ \phi_\ell(x) - \ba r_p) \, d\mu(x) \\
&\qquad = \frac{1}{n^4} \sum_{i=1}^n \sum_{k=1}^n \sum_{\ell=1}^n \sum_{p=1}^n \Cov(Y_{i,k},Y_{\ell,p}).\qedhere
\end{align*}

\end{enumerate}
\end{proof}

\begin{customprop}{7.2} \label{concordance-bounds}
Assume that $\mu_1,\dots,\mu_n$ are not all point masses (so the denominators of $W_{\text{scale}}$, $W_{\text{ratings}}$ are nonzero).
\begin{enumerate}[label = (\alph*)]

\item
$$0 \leq W_{\text{scale}} \leq 1.$$
Equality on the left never occurs. Equality on the right occurs iff all $\mu_i$ are the same, up to translation by a constant.

\item
$$0 \leq W_{\text{ratings}} \leq 1.$$
Equality on the left occurs iff user opinions on every item are equally balanced, i.e.\ $R_0$ is constant. Equality on the right occurs iff all user ratings are the same, up to translation by a constant, i.e.\ $\phi_i = \phi_\ell + a_{i,\ell}$ for all $1 \leq i, \ell \leq M$.

\end{enumerate}
\end{customprop}

\begin{proof}
The non-negativity of variance implies that both quantities are $\geq 0$. The upper bounds both follow from applying Jensen's inequality:
\begin{enumerate}[label = (\alph*)]

\item $\Var(\wh \mu) = 0$ when $\wh \mu$ is a point mass. Since $F_{\wh \mu}^{-1} = \frac{1}{n} \sum_{i=1}^n F_i^{-1}$, this can only happen when all $\mu_i$ are point masses and they are all equal. For the upper bound, we calculate
\begin{align*}
\Var(\wh \mu) &= \int (F_{\wh \mu}^{-1}(y) - \ba r)^2 \, dy \\
&= \int \paren{\frac{1}{n} \sum_{i=1}^n (F_i^{-1}(y) - \ba r_i)}^2 \, dy \\
&\leq \int \frac{1}{n} \sum_{i=1}^n (F_i^{-1}(y) - \ba r_i)^2 \, dy \\
&= \frac{1}{n} \sum_{i=1}^n \Var(\mu_i),
\end{align*}
with equality if and only if $F_i^{-1}(y) - \ba r_i = F_k^{-1}(y) - \ba r_k$ for all $1 \leq i,k \leq n$ and Lebesgue-a.e.\ $y$. This is precisely when $\mu_1 , \ldots , \mu_n$ are translates of each other.

\item
\begin{align*}
\Var((R_0)_* \mu) &= \int \paren{\frac{1}{n} \sum_{i=1}^n F_{\wh \mu}^{-1} \circ F_i \circ \phi_i(x) - \ba r }^2 \, d\mu(x) \\
&= \int \paren{\frac{1}{n^2} \sum_{k=1}^n \sum_{i=1}^n F_k^{-1} \circ F_i \circ \phi_i(x) - \ba r_k }^2 \, d\mu(x) \\
&\leq \int \frac{1}{n^2} \sum_{k=1}^n \sum_{i=1}^n (F_k^{-1} \circ F_i \circ \phi_i(x) - \ba r_k )^2 \, d\mu(x) \\
&= \frac{1}{n^2} \sum_{k=1}^n \sum_{i=1}^n \Var(\mu_k) \\
&= \frac{1}{n} \sum_{k=1}^n \Var(\mu_k),
\end{align*}
with equality if and only if $F_k^{-1} \circ F_i \circ \phi_i(x) - \ba r_k = F_p^{-1} \circ F_\ell \circ \phi_\ell(x) - \ba r_p$ for all $1 \leq i,k,\ell,p \leq n$ and $\mu$-a.e.\ $x$. This implies that $F_k^{-1},F_p^{-1}$ give the same distribution, except perhaps translated by a constant. The translation equivalence of these distributions, along with these $\mu$-a.e.\ equalities then implies that $\phi_i,\phi_\ell$ must give the same ordering to all $M$ atoms of $\mu$. Therefore, equality occurs exactly when all user rating profiles are the same, up to translation by a constant.\qedhere

\end{enumerate}
\end{proof}

\section{Calculations for Example 6}

\begin{customex}{6}
Let $\mu$ be any probability measure on $[0,1]$ which is symmetric about $1/2$ (besides the point mass at $1/2$), and let $\phi_k(x) := \alpha_k(x-1/2) + 1/2$ for some iid $\alpha_k$ with $\E[|\alpha_1|] = 1$ and $\P(\alpha_1 > 0) > 1/2$. Then, recalling from Example~\ref{scaling_and_reversing} that
$$R_0(x) = \paren{\frac{1}{n} \sum_{j=1}^n |\alpha_k|} \paren{\frac{1}{n} \sum_{k=1}^n \sgn(\alpha_k)}(x-1/2) + 1/2,$$
$$F_{\wh \mu}^{-1} \circ F_\mu = \paren{\frac{1}{n} \sum_{j=1}^n |\alpha_k|}(x-1/2) + 1/2,$$
we can calculate that
$$\Var(\mu_k) = \alpha_k^2\Var(\mu),$$
$$\Var((R_0)_*\mu) = \paren{\frac{1}{n} \sum_{j=1}^n |\alpha_k|}^2 \paren{\frac{1}{n} \sum_{k=1}^n \sgn(\alpha_k)}^2 \Var(\mu),$$
$$\Var(\wh \mu) =  \paren{\frac{1}{n} \sum_{j=1}^n |\alpha_k|}^2 \Var(\mu).$$
So we get that
$$W_{\on{scale}} = \frac{(\frac{1}{n} \sum_{k=1}^n |\alpha_k|)^2}{\frac{1}{n} \sum_{k=1}^n \alpha_k^2} \xrightarrow{a.s., n \to \infty} \frac{1}{\E[\alpha_1^2]},$$
$$W_{\on{ratings}} = \frac{(\frac{1}{n} \sum_{k=1}^n |\alpha_k|)^2 (\frac{1}{n} \sum_{k=1}^n \sgn(\alpha_k))^2}{\frac{1}{n} \sum_{k=1}^n \alpha_k^2} \xrightarrow{a.s., n \to \infty} \frac{\E[\sgn(\alpha_1)]}{\E[\alpha_1^2]}.$$
Since the distribution of $\alpha_k$ is subject to the constraint that $\E[|\alpha_k|] = 1$, we can see that $W_{\on{scale}}$ decreases as $\Var(|\alpha_1|)$ increases. Similarly, $W_{\on{ratings}}$ decreases as $\Var(|\alpha_1|)$ increases but also incorporates a multiplicative factor which decreases the value as the population agrees less on preference, i.e. when $\P(\alpha_1 < 0)$ grows closer to $1/2$.
\end{customex}

\section{Proofs of Lemma 5.4 and Theorem 5.5} \label{GC-appendix}

Here, we provide the proofs of Lemma~\ref{inverse_CDF_GC} and Theorem~\ref{wass-gc}. The latter will be a consequence of the former.

\begin{customlem}{5.4}[Glivenko--Cantelli-type theorem for inverse CDFs of Wasserstein barycenters on ${[0,1]}$]
\label{inverse_CDF_GC}
Let $\Lambda$ be a random measure on $[0,1]$, and let $\lambda$ be the population $W_2$-barycenter of the law of $\Lambda$. Draw iid samples $\Lambda_1,\Lambda_2,\dots$ according to the law of $\Lambda$, and let $ \lambda_n$ be the empirical $W_2$-barycenter of $\Lambda_1,\dots,\Lambda_n$. Denote the cumulative distribution functions of $ \lambda_n$ and $\lambda$ as $F_{ \lambda_n}, F_\lambda$, respectively, and define the generalized inverse cumulative distribution functions as
$$F_{ \lambda_n}^{-1}(x) := \inf \{ t : F_{ \lambda_n}(t) \geq x\}, \qquad F_{\lambda}^{-1}(x) := \inf \{ t : F_{\lambda}(t) \geq x\}.$$
Then
$$\sup_{x \in [0,1]} | F_{ \lambda_n}^{-1}(x) - F_{\lambda}^{-1}(x)| \to 0$$
almost surely as $n \to \infty$. Moreover, if $\supp \lambda = [0,1]$, then
$$\E \squa{\sup_{x \in [0,1]} | F_{ \lambda_n}^{-1}(x) - F_{\lambda}^{-1}(x)|} \leq 4 \sqrt{\frac{\log n}{n}}.$$
for all $n \geq 3$.
\end{customlem}

The key ingredient is the following pointwise convergence lemma for inverse CDFs of Fr\'echet means.

\begin{lemma} \label{pointwise_quantile_convergence}
Let $\Lambda$ be a random measure in $\mc P_2(\R)$ with $\E[W_2(\Lambda,\mu)] < \infty$ for all $\mu \in \mc P_2(\R)$. Then the population Fr\'echet mean $\lambda$ of the law of $\Lambda$ satisfies $F_\lambda^{-1}(x) = \E[F_\Lambda^{-1}(x)]$ for all $x \in (0,1)$.
\end{lemma}

The proof can be found in Chapter 3 of \cite{panaretos2020invitation}. The idea is that, since we know that $\| F_{\lambda_n}^{-1} - F_\lambda^{-1} \|_{L^2(dx)} = W_2(\lambda_n, \lambda) \to 0$, we can upgrade the $L^2$ convergence into pointwise convergence to an everywhere well-defined function by leveraging the fact that these functions are all nondecreasing and left continuous.

The proof of the a.s.\ convergence in Lemma~\ref{inverse_CDF_GC} follows from an application of the following deterministic lemma, which is usually used to provide a simple analytic proof of the classical Glivenko--Cantelli theorem.

\begin{lemma} \label{deterministic-gc}
Let $G,G_1,G_2,\dots$ be non-decreasing functions $[0,1] \to [0,1]$, which are all right continuous (resp. left continuous). If
\begin{enumerate}[label = (\roman*)]

\item $G_n(x) \to G(x)$ at each rational $x \in [0,1]$,

\item $G_n(x^-) \to G(x^-)$ (resp. $G_n(x^+) \to G(x^+)$) for all discontinuity points $x \in [0,1]$ of $G$,

\end{enumerate}
then $\sup_{x \in [0,1]} |G_n(x) - G(x)| \to 0$.
\end{lemma}

We will not provide a full proof of this deterministic lemma because the proof of the quantitative bound we presently prove for Lemma~\ref{inverse_CDF_GC} uses essentially the same argument (with some more bookkeeping for the a.s.\ convergence if $G$ has any discontinuity points).

\begin{proof}[Proof of Lemma~\ref{inverse_CDF_GC}]
Since $\supp \lambda = [0,1]$, $F_\lambda^{-1}$ is continuous (and nondecreasing), $F_\lambda^{-1}(0) = 0$, and $F_\lambda^{-1}(1) = 1$. Thus, there exist points $0 = t_1 < \cdots < t_N = 1$ such that $F_\lambda^{-1}(t_{k+1}) - F_\lambda^{-1}(t_k) \leq \sqrt{\frac{\log n}{2n}}$ for all $1 \leq k \leq N-1$. Necessarily, $N \leq \sqrt{\frac{2n}{\log n}}$. By Lemma~\ref{pointwise_quantile_convergence} and Hoeffding's inequality for bounded random variables, for each $k$,
$$\P\paren{|F_{\lambda_n}^{-1}(t_k) - F_\lambda^{-1}(t_k)| > \sqrt{\tfrac{\log n}{2n}}} \leq 2\exp \paren{-2n \paren{\sqrt{\frac{\log n}{2n}}}^2} = \frac{2}{n}.$$
Then
$$\P \paren{|F_{\lambda_n}^{-1}(t_k) - F_\lambda^{-1}(t_k)| > \sqrt{\tfrac{\log n}{2n}} \text{ for any $k$}} \leq N \frac{2}{n} \leq \frac{2^{3/2}}{\sqrt{n\log n}}.$$

We now show that if $|F_{\lambda_n}^{-1}(t_k) - F_\lambda^{-1}(t_k)| \leq \sqrt{\frac{\log n}{2n}}$ for all $k$, then we have $\sup_{x \in [0,1]} |F_{\lambda_n}^{-1}(x) - F_\lambda^{-1}(x)| \leq \sqrt{\frac{2\log n}{n}}$,
so that
$$\P \paren{\sup_{x \in [0,1]} |F_{\lambda_n}^{-1}(x) - F_\lambda^{-1}(x)| > \sqrt{\tfrac{2\log n}{n}}} \leq \frac{2^{3/2}}{\sqrt{n\log n}}.$$
Indeed, for a point $x$ with $t_{j-1} < x < t_j$, given this event for all $j$,
\begin{align*}
&|F_{\lambda_n}^{-1}(x) - F_\lambda^{-1}(x)| \\
&\quad \leq \begin{cases}
F_{\lambda_n}^{-1}(t_j) - F_\lambda^{-1}(t_{j-1}) &\text{if } F_{\lambda_n}^{-1}(x) \geq F_\lambda^{-1}(x) \\
F_\lambda^{-1}(t_j) - F_{\lambda_n}^{-1}(t_{j-1}) &\text{if } F_{\lambda_n}^{-1}(x) < F_\lambda^{-1}(x)
\end{cases} \\
&\quad \leq \begin{cases}
|F_{\lambda_n}^{-1}(t_j) - F_\lambda^{-1}(t_j)| + |F_\lambda^{-1}(t_j) - F_{\lambda}^{-1}(t_{j-1})| &\text{if } F_{\lambda_n}^{-1}(x) \geq F_\lambda^{-1}(x) \\
|F_\lambda^{-1}(t_j) - F_\lambda^{-1}(t_{j-1})| + | F_\lambda^{-1} (t_{j-1}) - F_{\lambda_n}^{-1}(t_{j-1})| &\text{if } F_{\lambda_n}^{-1}(x) < F_{\lambda}^{-1}(x)
\end{cases} \\
&\quad \leq 2 \sqrt{\frac{\log n}{2n}}.
\end{align*}
We now get the $L^1$ bound by
\begin{align*}
& \E \squa{\sup_{x \in [0,1]} |F_{\lambda_n}^{-1}(x) - F_\lambda^{-1}(x)| } \\
& \qquad \leq \sqrt{\frac{2\log n}{n}} \cdot \P \paren{\sup_{x \in [0,1]} |F_{\lambda_n}^{-1}(x) - F_\lambda^{-1}(x)| \leq \sqrt{\tfrac{2\log n}{n}}} \\
& \qquad \qquad + 1 \cdot \P \paren{\sup_{x \in [0,1]} |F_{\lambda_n}^{-1}(x) - F_\lambda^{-1}(x)| > \sqrt{\tfrac{2\log n}{n}}} \\
&\qquad \leq \sqrt{\frac{2\log n}{n}} + \frac{2^{3/2}}{\sqrt{ n \log n}} \\
&\qquad \leq \paren{\sqrt 2 + \frac{2^{3/2}}{\log 3}} \sqrt{\frac{\log n}{n}} \\
&\qquad \leq 4 \sqrt{\frac{\log n}{n}}
\end{align*}
where the second to last inequality holds for all $n \geq 3$.
\end{proof}

The constant factor in the bound can be decreased to any value $ > \sqrt 2$, at the cost of increasing the minimum value of $n$ to which the inequality applies.

We can, without much extra work, leverage this uniform convergence of inverse CDFs to obtain uniform convergence of the actual CDFs. The following result will not be used in this paper but is interesting in its own right.

\begin{customthm}{5.5}[Glivenko--Cantelli-type theorem for Wasserstein barycenters on ${[0,1]}$] \label{wass-gc}
Let $\Lambda$ be a random measure on $[0,1]$, and let $\lambda$ be the population $W_2$-barycenter of the law of $\Lambda$. Draw iid samples $\Lambda_1,\Lambda_2,\dots$ according to the law of $\Lambda$, and let $ \lambda_n$ be the empirical $W_2$-barycenter of $\Lambda_1,\dots,\Lambda_n$. Denote the cumulative distribution functions of $ \lambda_n$ and $\lambda$ as $F_{ \lambda_n}, F_\lambda$, respectively. Assume that $\lambda$ has no atoms and that $\supp \Lambda = [0,1]$ a.s. Then
$$\sup_{x \in [0,1]} | F_{ \lambda_n}(x) - F_{\lambda}(x)| \to 0$$
almost surely as $n \to \infty$.
\end{customthm}

We will use our previous results, along with the following deterministic lemma.

\begin{lemma} \label{deterministic-inverse}
Let $f_1,f_2,\ldots : [0,1] \to [0,1]$ be a sequence of functions with right inverses $f_n^{-1}$ (i.e.\ $f_n \circ f_n^{-1} = \id$) converging uniformly to a function $f : [0,1] \to [0,1]$ with continuous left inverse $f^{-1}$ (i.e.\ $f^{-1} \circ f = \id$). Then $f_n^{-1} \to f^{-1}$ uniformly, as well.
\end{lemma}

\begin{proof}[Proof of Lemma~\ref{deterministic-inverse}]
Since $f^{-1}$ is continuous on a compact domain, it is uniformly continuous. So, given $\varepsilon > 0$, let $\delta > 0$ be such that $|x-y| < \delta \implies |f^{-1}(x) - f^{-1}(y)| < \varepsilon$ for all $x,y \in [0,1]$. Then if $n$ is large enough such that $\sup_{y \in [0,1]} |f_n(y) - f(y)| < \delta$,
\begin{align*}
|f_n^{-1}(x) - f^{-1}(x)| &= |f^{-1}(f(f_n^{-1}(x)) - f^{-1}(x)| \\
&= |f^{-1}(f(f_n^{-1}(x)) - f^{-1}(f_n(f_n^{-1}(x)))|.
\end{align*}
$|f(f_n^{-1}(x)) - f_n(f_n^{-1}(x))| < \delta$, so the above is $< \varepsilon$.
\end{proof}

\begin{proof}[Proof of Theorem~\ref{wass-gc}]
We will apply Lemma~\ref{deterministic-inverse} to $f_n = F_{\lambda_n}^{-1}$ and $f = F_\lambda^{-1}$ (verifying the uniform convergence of the $F_{\lambda_n}^{-1}$ via Lemma~\ref{inverse_CDF_GC}) to get that $F_{\lambda_n} \to F_\lambda$ uniformly $\P$-a.s. We now verify the conditions of the lemma.

First, the atomless condition on $\lambda$ guarantees that $F_\lambda$ is continuous. It is also true by definition that $F_\lambda(F_\lambda^{-1}(x)) = x$. In general, $F_{\lambda_n}^{-1}(F_{\lambda_n}(x)) = \inf \{ y : F_{\lambda_n}(y) = F_{\lambda_n}(x) \} \leq x$, but we have equality if $F_{\lambda_n}$ is strictly increasing, i.e.\ $\supp \lambda_n = [0,1]$. This follows from the fact that $\supp \Lambda = 1$ a.s., as in general,
$$F_\mu \text{ is strictly increasing} \iff \supp \mu = [0,1] \iff F_{\mu}^{-1} \text{ is continuous}.$$
Averaging continuous functions yields a continuous function. So $F_{\lambda_n}^{-1} := \frac{1}{n} \sum_{k=1}^n F_{\Lambda_k}^{-1}$ is continuous, meaning $F_{\lambda_n}$ is strictly increasing.
\end{proof}
\end{appendix}

\pagebreak

\bibliographystyle{alpha}
\bibliography{rating-scales}

\end{document}